\documentclass[final,3p]{elsarticle}

\usepackage[colorlinks=true,breaklinks=true,pdftex]{hyperref}

\usepackage{amsthm}
\usepackage{amsfonts,amsmath,amssymb,amstext}
\usepackage{graphicx}
\usepackage{latexsym}
\usepackage{setspace}

\usepackage{soul}
\renewcommand{\hl}[1]{#1}

\theoremstyle{plain}% default

\theoremstyle{definition}
\newtheorem{definition}{Definition}[section]
\newtheorem{example}{Example}[section]

\theoremstyle{remark}

\def\Re{\mathbb{R}}

\def\R{\mathsf{R}}
\def\W{\mathsf{W}}
\def\w{\mathsf{w}}
\def\fv{\mathsf{fv}}
\def\A{\mathsf{A}}
\def\inn{\iota}
\def\bd{\beta}
\def\diag{\mathsf{diag}}
\def\diags{\mathsf{diags}}
\newcommand{\norm}[1]{\left\lVert#1\right\rVert}

\journal{arXiv}

\biboptions{authoryear}

\begin{document}

\begin{frontmatter}

\title{Domain Decomposition for Mean Curvature Flow of Surface Polygonal Meshes}

% USE FOR THE REVIEW PROCESS
% \author{1076}

%% ONLY USE FOR THE FINAL ACCEPTED PAPER SUBMISSION
\author[first]{Lenka Ptáčková\corref{cor}}
\ead{lenka@kam.mff.cuni.cz}
\cortext[cor]{Corresponding author}
\author[first]{Michal Outrata}
\ead{outrata@karlin.mff.cuni.cz}
\address[first]{Department of Numerical Mathematics, Charles University, Sokolovsk\'a 49/83, Prague, 186 75, Czech Republic}

\begin{abstract}
We examine the use of domain decomposition for potentially more efficient mean curvature flow of surface meshes, whose faces are arbitrary simple polygons. We first test traditional domain decomposition methods with and without overlap of deconstructed domains. And we present adapted Robin transmission conditions of optimized Schwarz method. We then analyze the resulting smoothing from the point of view of shape quality and texture deformation. By decomposing the initial mesh into two sub-meshes, we solve two smaller boundary value problems instead of one big problem, and we can process these two tasks almost entirely in parallel.
\end{abstract}


\end{frontmatter}

% Comment out for final accepted paper submission
% \linenumbers

%%%%%%%%%%%%%%%%%%%%%%%%%%%%%%%%%%%%%%%%%%%%%%%%%%%%%%%%%%%%%%%%%%%%%
%%%%%%%%%%%%%%%%%%%%%%%%%%%%%%%%%%%%%%%%%%%%%%%%%%%%%%%%%%%%%%%%%%%%%
\begin{figure}[ht!]
\centering
\includegraphics[width=0.19\textwidth]
{ 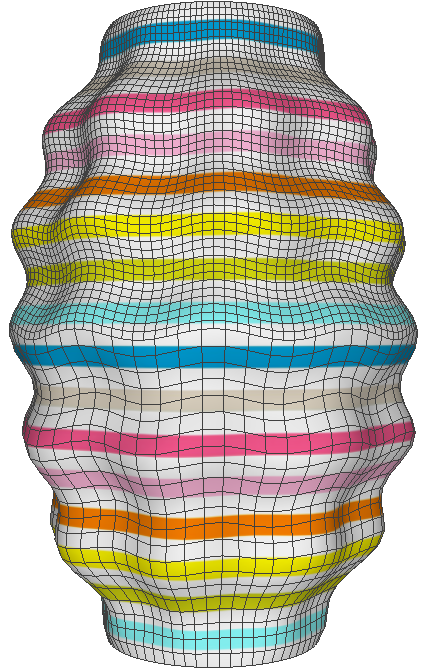}
\includegraphics[width=0.19\textwidth]
{ 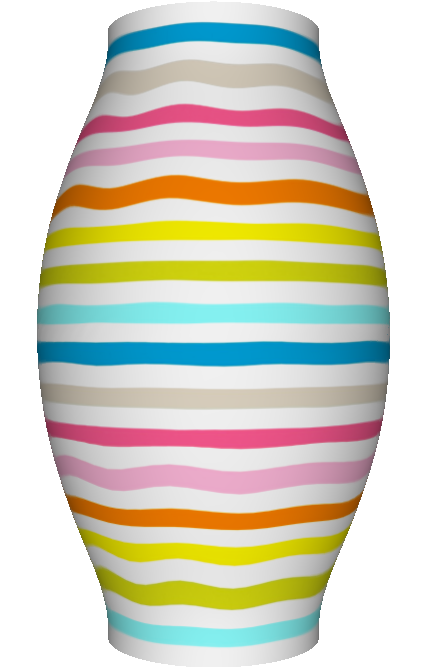}
\includegraphics[width=0.19\textwidth]
{ 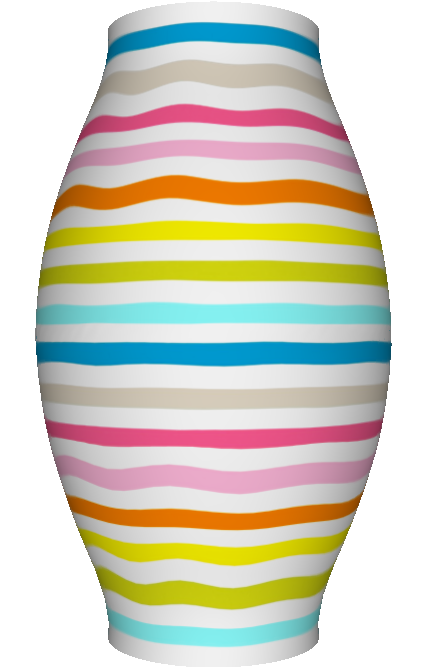}
\includegraphics[width=0.19\textwidth]
{ 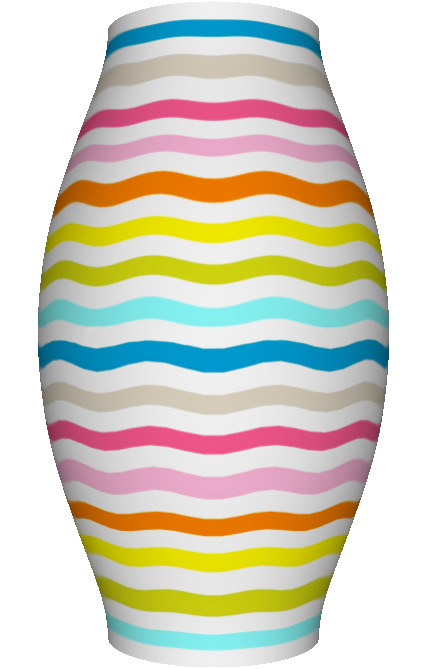}
\includegraphics[width=0.19\textwidth]
{ 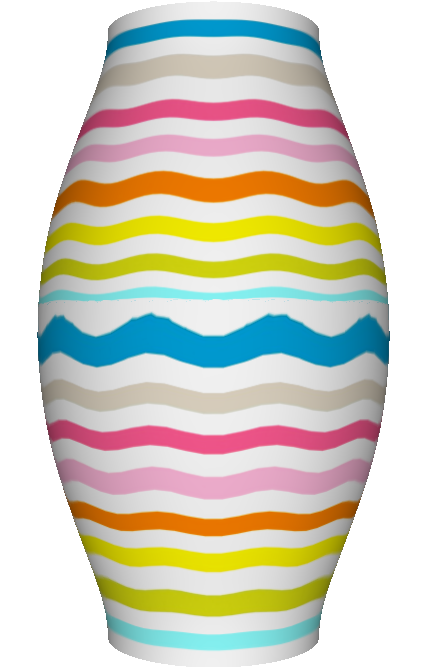}
\centering
\caption{Mesh smoothing through mean curvature flow (MCF) of a surface polygonal mesh (far left), whose upper part has the double of the level of detail of the lower part. We decompose the mesh into two sub-meshes (separate the upper and the lower part). We then apply five iterations of MCF with our adapted Robin transmission conditions and time step $dt=0.05$. The show the resulting meshes for Laplacians of \cite{AlexaWardetzky2011} (second image), \cite{PtackovaVelho2021} (third image), and \cite{Fujiwara1995} (fourth image). The last fifth image is the resulting mesh after MCF without decomposition for the Laplacian of \cite{Fujiwara1995}.
We can see that even though the level of details of the sub-meshes is different and the interface is a curved line, our adapted Robin transmission conditions work well and the resulting meshes are well smoothed. In the case of the Laplacian of \cite{Fujiwara1995}, our domain decomposition method actually prevents the tangential shifting of vertices caused by the different levels of details, which can be observed as the significant texture deformation in the fifth image.}
\label{fig:teaser}
\end{figure}

\section{Introduction}

Domain decomposition methods (DDM) have been traditionally used for solving boundary value problems over a given domain by splitting the domain into smaller subdomains and subsequently solving smaller problems on these subdomains. This technique has a long history, the German mathematician Hermann Schwarz is considered to be the father of the idea. In 1869-70 he introduced the Schwarz alternating method as a tool for solving Laplace equation on complex domains. Nowadays, in the age of parallel computing, these methods are becoming more and more pertinent. In this paper, we examine the application of DDM in an unconventional context of mean curvature flow of surface polygonal meshes.

At its core, mean curvature flow (MCF) of surface meshes with a boundary corresponds to solving the heat equation with specified boundary conditions. We formalize the statement of our problem in Section \ref{sec:MCF}.
On differential manifolds, MCF of surfaces would eventually result in minimal surfaces with given (fixed) boundary. Finding such a minimal surface corresponds to solving the Dirichlet problem. However, we focus on MCF and its (temporary) smoothing effect, that is, we are interested in the MCF for a limited time frame.

Schwarz and optimized Schwarz DDM, concretely with Robin transmission conditions, have been traditionally used in solving the Poisson equation with boundary conditions. However, since we are interested in the process of mean curvature flow itself, we must adapt the transmission conditions --- this adaptation and its analysis are our main contributions. We present them in Section \ref{sec:mcfDDM}.
In order to employ \hl{these} conditions not only grid--like meshes, but also on general polygonal meshes, we present a novel normal derivative in Definition \ref{def:newNormalDerivative}.

Furthermore, we employ various discrete Laplace operators in order to prove that our method is not tailored for one specific discretization of the Laplacian. We compare the resulting flow from various perspectives. The most straightforward viewpoint is the quality of resulting smoothing, i.e., minimization of negative artifacts introduced by domain decomposition. However we also analyze the influence of the decomposition on the deformation of a texture mapped on the smoothed surface --- compared to the smoothing without domain decomposition.

We evaluate the adequacy of the proposed method through series of numerical tests. In each test, the reference solution is the mesh smoothed without DDM, and two sets of measures are used: the pseudo--Hausdorff distance of the corresponding meshes (which quantifies the quality of overall smoothing) and the relative distance of corresponding vertices (which quantifies the deformation introduced by the decomposition).

The proposed adaptation seems to be very natural and yet, to the best of our knowledge, it is new. Indeed, DDM have been extensively investigated previously, but not in the setting of mesh smoothing. Similarly, mesh smoothing by mean curvature flow is an important and well studied topic in the computer graphics community. Before we present our framework, we provide a review of related literature in the following section.

%===============================================================

\subsection{Related literature}\label{sec:RelatedLiterature}
DDM have been extensively used for solving Poisson's equation with various boundary conditions --- finding a minimal surface with a given boundary is an example of such a problem.

\subsubsection*{Minimal surfaces}
A minimal surface is a surface whose mean curvature is equal to zero at any interior point. Showing existence of a minimal surface with a given closed boundary curve is also known as the Problem of Plateau. Its name honors the work of the Belgian physicist J. A. Plateau, who performed extensive experimental studies with soap bubbles.

One of the earliest and undoubtedly very influential works addressing the problem in discrete setting is the seminal paper of\hl{ }\cite{PinkallPolthier93}. They assume a Dirichlet problem, i.e., given a triangle mesh with boundary, find a new mesh such that the Laplacian of vertex position at each interior vertex is equal to zero, while the boundary points stay fixed or are allowed to move tangentially along the boundary. They employed what was later coined as the cotan--Laplace formula for discrete Laplacian on 0--forms.

\subsubsection*{Discrete mean curvature flow}
Mean curvature flow of discrete surfaces, i.e., surface meshes, in discrete time has been traditionally used for mesh smoothing. As such, MCF is often employed in order to remove rough features coming from irregularly triangulated data. For example, we might want to remove some undesirable noise (high frequencies of the mesh), while maintaining the desirable geometric features (low frequencies) of the initial mesh.
We adopt the approach of \cite{Desbrun1999} in the sense that we employ a discrete Laplace operator together with the backward Euler method to find the new vertex positions of the smoothed mesh. We extend their scheme from purely triangular meshes to any polygonal mesh --- and introduce domain decomposition methods.

\subsubsection*{Domain decomposition methods}
Given a partial differential equation (PDE), or a system of linear algebraic equations, most DDM decompose the domain into several subdomains, or the set of unknowns into several subsets, and formulate new subdomain problems posed only on these subdomains, formulated separately for the unknown subsets. As the resulting subdomain solutions usually do not coincide with the restrictions of the solution of the original problem, we introduce an iterative process of \emph{information exchange}. Thus we solve the subdomain problems repeatedly while exchanging information between the subdomains. Choosing an apt information exchange scheme, the error of the subdomain solutions decreases with each iteration and the subdomain solutions converge to the restrictions of the original one, see~\cite{Dolean2015Intro, Gander2024TimeParallel} for further details.

Notably, the so-called \textbf{(optimized) Schwarz waveform relaxation methods}\hl{,} see~\cite[Chapter 3]{Gander2024TimeParallel}\hl{,} apply this idea to the heat equation, which is of interest in our setting. There the information is exchanged by overlapping subdomains and/or by appropriate interface conditions. It is well known, within the DDM community, that the larger the overlap, the faster is the convergence; but also the higher are the computational expenses. Unfortunately, in our setting a simple exchange of information by overlap does not give satisfactory results if the interface is placed on a peak. Since even for large overlaps, if the interface is placed in an area with high mean curvature, the surface is not smoothed evenly.

Thus we must \hl{resort} to more complex information to be exchanged: standard conditions are combined \emph{Dirichlet} and \emph{Neumann}, corresponding to \emph{Robin}, as well as \emph{Ventcell} transmission conditions.
They enforce compatibility of a linear combination of the solution at the interface, the normal derivative along the interface, and the tangential second-order derivative along the interface of the subdomain solutions. The choice of the coefficients of this linear combination then plays a decisive role in the quality of the resulting method and these are, naturally, problem dependent.

Moreover, the use of these methods for MCF poses new, intriguing challenges. The DDM global approximations are in essence non--smooth or even discontinuous due to the imperfect information exchange --- only as the method converges to the solution we get smoother approximations. However, since we are interested in the smoothing flow itself, we must adapt the transmission conditions to guarantee the continuity and the desired degree of smoothness after each iteration of MCF. We present these modifications in Section \ref{sec:mcfDDM}.

\subsubsection*{Domain decomposition methods in computer graphics}
The domain decomposition methods are not new to the Computer Graphics (CG) community and they have been successfully applied in tasks ranging from fluid flow \citep{English2013} to elastic simulations \citep{MGLNF2015}.

Multigrid methods are another example of the use of DDM in CG, see \citep{Wiersma2023, Liu2021}. For instance, authors of \citep{ChuangKazhdan2011} use a finite-elements hierarchy that supports a multigrid solver for performing the semi-implicit MCF.

However all of these geometry-based approaches \citep{MGLNF2015, Wiersma2023, Liu2021, English2013, ChuangKazhdan2011} use an additional structure such as the multigrid or the Chimera grids. Whereas we use only the given mesh and employ the standard Schwarz and our adapted optimized Schwarz method.

Our approach to DDM in MCF is thus much more straightforward and is similar to the Boundary-Only Coupling of \cite{sellan2019solid}, which can be seen as a standard Schwarz method, i.e., not optimized, and instead become interested in the problem of numerous subdomains; while working with tetrahedral meshes.

Similarly, in~\cite{zhou2024alternating} the authors also consider a ``classical'' variant of the domain decomposition, in the sense of exchanging only the Dirichlet trace data. Many other works, including~\cite{benek1983flexible,henshaw1994fourth,dobashi2008fast}, consider splitting computational domains or domains of discretization into multiple subdomains and consider the resulting discretized coupling.

One of the novelties we bring is testing the results about Schwarz methods in the area of mesh smoothing, (numerically) exploring to what extent can these be relied upon in this context.

Even though our adaptation (Section \ref{subsec:OurAlgorithm}) looks very natural, we are not aware of any publication presenting it. Moreover, to the best of our knowledge, there is no publication illustrating the direct use of standard and optimized Schwarz methods directly to the MCF of general polygonal meshes in such a compact manner.

Let us summarize the key differences between our approach and the methods presented in the previous CG literature or finite elements literature:
\begin{enumerate}
\item We work with general polygonal meshes, i.e., our framework is not limited to pure triangle or pure quad meshes.
\item We do not involve any additional structure or mesh as in multigrid methods or Chimera grids over a curved surface mesh.
\item Our adaptation (Section \ref{subsec:OurAlgorithm}) does not involve any overlap and thus it does not introduce unnecessary degrees-of-freedom. And it can even give better results than MCF without domain decomposition, as illustrated in Figure \ref{fig:teaser}.
\end{enumerate}

%===============================================================
%===============================================================

\section{Preliminaries}\label{sec:MCF}

Let $\{M_t \}_{t\in\Re}$ be an evolving family of smooth surfaces with or without boundary in $\Re^3$. We say that $M_t$ is moving by mean curvature in time $t$ if it satisfies the nonlinear parabolic equation
\[
 \frac{\partial}{\partial t}x = \vec{H}(x),\;\; x\in M_t,\;\; t\in\Re,
\]
where $\vec{H}(x)$ is the mean curvature vector of $M_t$ at point $x$. The mean curvature vector is defined by
\[
 \vec{H} = H\vec{n} = (\kappa_1 + \kappa_2)\vec{n},
\]
where $\kappa_1,\kappa_2$ are the principal curvatures, $H$ is the mean curvature, and $\vec{n}$ is the unit normal of the given surface.
However, at an interior point on the surface, the arithmetic mean of the normal curvatures in any two orthogonal directions is equal to the mean curvature at that point \citep{Araujo2024}.

\subsection{Mean curvature flow as a heat equation in discrete setting}
Let us denote $\Delta$ the Laplace--de Rham operator on functions $f : M\to\Re$ on manifold $M\subset \Re^3$. Then for local coordinates $X$ of points on $M$ we have
\begin{equation}
 \Delta X (x) = \vec{H}(x),\quad x \in M\setminus \partial M.
\end{equation}

The mean curvature flow, which aims at decreasing the magnitude of the mean curvature, can be then written as a \textbf{heat equation}:
\begin{equation}\label{eq:heatEquation}
 \frac{\partial X}{\partial t} = - \Delta X.
\end{equation}

In our discrete setting, an orientable surface is represented as a polygonal surface mesh, or equivalently, as an orientable 2--dimensional pseudomanifold, whose faces are any simple polygons, possibly non--planar and non--convex. Thus $X$ corresponds to coordinates of vertices of our mesh. Moreover, we work with discrete time, the explicit (forward) Euler scheme thus reads:
\begin{equation*}
 V^{n+1} = (I - dt \Delta) V^n,
\end{equation*}
with $V^0$ corresponding to the vertex coordinates of the initial mesh $M$ and $dt$ denoting the discrete time step.

However, the stability criterion of explicit Euler scheme requires smaller steps. Concretely, authors of \cite{Desbrun1999} assert that the stability criterion requires $dt \leq min(\norm{e}^2)$, where $\norm{e}$ is the edge length. Thus to produce noticeable smoothing or to converge to minimal surface with given boundary curves, a large number of iterations must be done, and in each iteration, the Laplace operator must be computed --- which demands integration of mesh elements. This is computationally expensive.

An implicit time integration of the heat equation (\ref{eq:heatEquation}) results in a numerically much more stable curvature flow. Now the new positions are computed as $V^{n+1} = V^n - dt \Delta V^{n+1}$. In this \textbf{backward Euler method} we find the updated vertex positions $V^{n+1}$ as solutions of the following linear system:
\begin{equation}\label{eq:MCFbackwardEuler}
(I + dt \Delta) V^{n+1} = V^n.
\end{equation}
In the above equation, the operator $\Delta$ should be computed based on the vertex coordinates $V^{n+1}$. However, in practice it is approximated by $\Delta$ computed relative to the known positions $V^n$, although the values of matrices $\Delta$ change after each iteration, since the vertex positions change.

%===============================================================

\subsection{Discrete Laplacians}
In order to show that our adaptation of the optimized Schwarz domain decomposition method is not tailored for one specific discretization of the Laplace operator, we work with three different discretizations of the Laplacian; all of which are defined on general polygonal meshes. We shall present them in matrix form. We thus start with definition of a set of \textit{building matrices}, \hl{i.e., matrices whose multiplication leads} to the discrete Laplacians\hl{. T}hese \textit{building matrices} are later employed also in computing normal derivatives in Section \ref{subsec:OptimizedSchwarz}.

Let $M$ be a mesh with vertices $V$, oriented halfedges $E$, and oriented faces $F$. Let also denote $\prec, \succ$ the incidence relation, i.e., $f\succ v$ means that $v$ is a vertex of face $f$. Further, incidence number $[f_i:e_j] = 1$ means that $e_j$ is a halfedge incident to $f_i$ with the same orientation as $f_i$. Since we work with halfedge structure, we have that either $[f_i:e_j] = 1$ or $[f_i:e_j] = 0$.

We denote $d_0\in\Re^{|E|\times|V|}$ and $d_1\in\Re^{|F|\times|E|}$ matrices representing discrete exterior derivatives on 0--forms and 1--forms, resp.,
\begin{equation}
d_0[i,j] = [e_i:v_j] =
\begin{cases}
1   & \text {if $v_j$ is the endpoint of $e_i$}, \\
-1  & \text {if $v_j$ is the starting point of $e_i$}, \\
0   & \text{otherwise,}
\end{cases}
\label{eq:d0}
\end{equation}

\begin{equation}
d_1[i,j] = [f_i:e_j] =
\begin{cases}
1   & \text{if $e_j\prec f_i$, $e_j$ is oriented as $f_i$}, \\
0   & \text{otherwise.}
\end{cases}
\label{eq:d1}
\end{equation}
Further, matrix $\fv\in\Re^{|F|\times|V|}$ expressing the incidence relation between faces and vertices reads
\begin{equation}
\fv[i,j] =
\begin{cases}
\frac{1}{p_i} & \text {if } v_j \prec f_i,\; f_i\text{ is a } p_i\text{--gon}, \\
0 & \text{otherwise.}
\end{cases}
\label{eq:fv}
\end{equation}
Matrix $\R\in\Re^{|E|\times |E|}$ then encodes the adjacency of edges, it is a block diagonal matrix given per $p_i$--gons $f_i$ as
\begin{eqnarray}
\R    &=& \diags(\R_{f_0},\dots,\R_{f_{|F|-1}})\\ \label{eq:R_start}
\R_f  &=&
\sum_{a=1}^{\lfloor\frac{p-1}{2} \rfloor} \bigg( \frac{1}{2} - \frac{a}{p} \bigg) \R_a,\\
\R_a[k,j] &=&
\begin{cases}
1  & \text{\small if $e_j$ is $(k+a)$--th halfedge of $f, [f:e_j]=1$}, \\
-1 & \text{\small if $e_j$ is $(k-a)$--th halfedge of $f, [f:e_j]=1$}, \\
0 & \text{\small otherwise. }
\end{cases}
\label{eq:R_end}
\end{eqnarray}
Matrix $\A\in\Re^{|E|\times |E|}$ then encodes whether an edge has two halfedges (it is an interior edge) or it has only one halfedge (it is a boundary edge). It serves for linearly combining the values of halfedges incident to the same edge:
\begin{equation}
\A[i,j]=
\begin{cases}
1  & \text{\small if $i=j$ and $e_j$ is a boundary edge}, \\
\frac{1}{2} & \text{\small if $i=j$ and $e_j$ is an interior edge}, \\
-\frac{1}{2} & \text{\small if $e_j$ is a halfedge of the same edge as $e_i$ but with orientation opposite to $e_i$}, \\
0 & \text{\small otherwise. }
\end{cases}
\label{eq:A}
\end{equation}

We also define metric--dependent matrices $\W_V\in\Re^{|V|\times|V|}$,
$\W_1\in\Re^{|E|\times|E|}$, and vector $\w_E\in\Re^{|E|\times 1}$ as:
\begin{equation}
\W_V[i,i] =\frac{1}{\sum\limits_{f_k\succ v_i}\frac{|f_k|}{p_k}},
\label{eq:WV}
\end{equation}
\begin{equation}
\W_1[i,j] =
\begin{cases}
\frac{\langle e_i,e_j \rangle}{|f|}
& \text{if $\exists f$ s.t. $[f:e_i] = 1$ and $[f:e_j]=1$}, \\
0   & \text{otherwise,}
\end{cases}
\end{equation}
\begin{equation}
\w_E = \Big(\frac{1}{\norm{e_0}_2},\dots,\frac{1}{\norm{e_{|E|-1}}}\Big)^\top, \text{ where }
(e_0,\dots,e_{|E|-1})^\top= d_0 V.
\end{equation}
Since we work with meshes, whose faces are $p$--gons, possibly non--planar and non--convex, we compute the area of a face $f$ as the magnitude of its vector area:
\[
|f| =\frac{1}{2} \norm{ \sum_{i=0}^{p-1} v_i\times v_{(i+1)\text{ mod }p}}_2 \; \text{ for } f= (v_0,\dots, v_{p-1}).
\]
The magnitude $|f|$ of the vector area of face $f$ is the largest area over all orthogonal projections of $f$ to planes in $\Re^3$. The vector area depends only on the boundary curve, not on the spanning surface \cite[Lemma 1]{AlexaWardetzky2011}.

%===============================================================

\textbf{The scale--dependent umbrella operator of}\hl{ }\cite{Fujiwara1995} is given by the following formula
\[
\Delta (v_i) = \frac{1}{m(v_i)}
\sum_{v_j\sim v_i} \frac{v_i - v_j}{\norm{v_i-v_j}_2},\quad
m(v_i) = \sum_{v_j\sim v_i} \norm{v_i-v_j}_2,
\]
where $\sum_{v_j\sim v_i}$ means to take the sum over all vertices $v_j\in V$ connected to $v_i$ by an edge. We adapt the Laplacian by multiplying it by the valence of $v_i$, so that it matches the standard second order five point finite difference discretization \cite[Section 2]{Gander2022} on uniform grids, up to a sign. Thus, we arrive at the following matrix form of the adapted Fujiwara's Laplacian
\begin{equation}
\Delta_F = \frac{1}{2}
\frac{|d_0^\top|\cdot (1,\dots,1)^\top}{|d_0^\top|\cdot \w_E}
d_0^\top\; \diag(w_E)\;\; d_0,
\label{eq:umbrellaLaplace}
\end{equation}
where the absolute value $|d_0^\top|$ is meant entry-wise.

\textbf{The cotan Laplacian and the polygonal Laplacian of}\hl{ }\cite{AlexaWardetzky2011}.
The cotan Laplacian \cite[eq. (6)]{PinkallPolthier93} is defined on a triangle mesh and takes into account not only the lengths of incident edges, but also the angles adjacent to $v_i$. Authors of \cite{AlexaWardetzky2011} present a class of discrete Laplacians on general polygonal meshes that leads to the cotan Laplacian on triangle meshes ---see Theorem 2 therein. We focus on their \textit{geometric} Laplacian, denoted here as $\Delta_A$, which can be expressed in matrix form using the previously given matrices as
\begin{equation}
\Delta_A = \W_V \; d_0^\top \;\R \W_1\R^\top d_0.
\label{eq:AWLaplace}
\end{equation}

\textbf{The polygonal Laplacian of}\hl{ }\cite{PtackovaVelho2021}
is actually a discrete Laplace--de Rham operator, also called the Hodge Laplacian, since it is defined using their discrete Hodge star operator. We denote their Laplacian $\Delta_P$ and it reads
\begin{equation}
 \Delta_P = -\W_V \fv^\top d_1 \A \W_1 \R^\top d_0.
\label{eq:PVLaplace}
\end{equation}
%===============================================================
%
Note that the all the above Laplacians correspond to the Laplace--de Rham operator, thus they are the negative of the Laplace--Beltrami operator, i.e., for a function $g$
\[
 \Delta g = -div(\nabla g).
\]

\subsection*{Boundary conditions}
If the underlying surface has a non--empty boundary $\partial M$, we keep the boundary vertices fixed. Denote `$\inn$' the indices of the interior vertices, i.e., $V_\inn$ is the set of interior vertices, and similarly $V_\bd$ is the set of all boundary vertices of $M$. For simplicity, let us assume that vertices of the mesh are ordered such that $V = (V_\bd,V_\inn)^\top$. The given Dirichlet boundary conditions are then represented as an identity sub-matrix, the equation (\ref{eq:MCFbackwardEuler}) then takes the form
\[
\begin{pmatrix}
I & 0\\
dt\Delta_{\bd} & I +dt\Delta_\inn
\end{pmatrix}
\begin{pmatrix}
V_\bd^n\\
V_\inn^{n+1}
\end{pmatrix}
=
\begin{pmatrix}
V_\bd^n\\
V_\inn^{n}
\end{pmatrix}.
\]
The Laplacian can be divided column--wise into two sub-matrices $\Delta = (\Delta_{\bd}|\Delta_\inn)$. We thus solve a smaller system reading
\begin{equation}\label{eq:mcfWithBoundary}
(I + dt\Delta_\inn)V_\inn^{n+1} = V_\inn^{n} - dt\Delta_\bd V_\bd^n.
\end{equation}
We sometimes omit the subscripts of $\Delta_\inn, \Delta_\bd$ to avoid cluttering of the notation.

\begin{figure}
\includegraphics[width=0.24\linewidth]{ 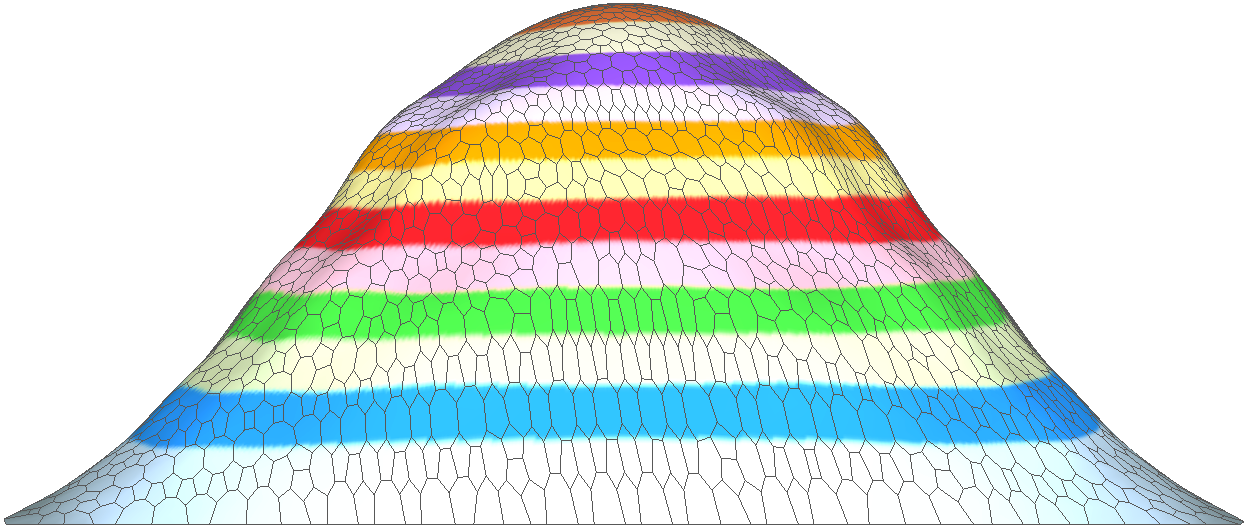}
\includegraphics[width=0.24\linewidth]{ 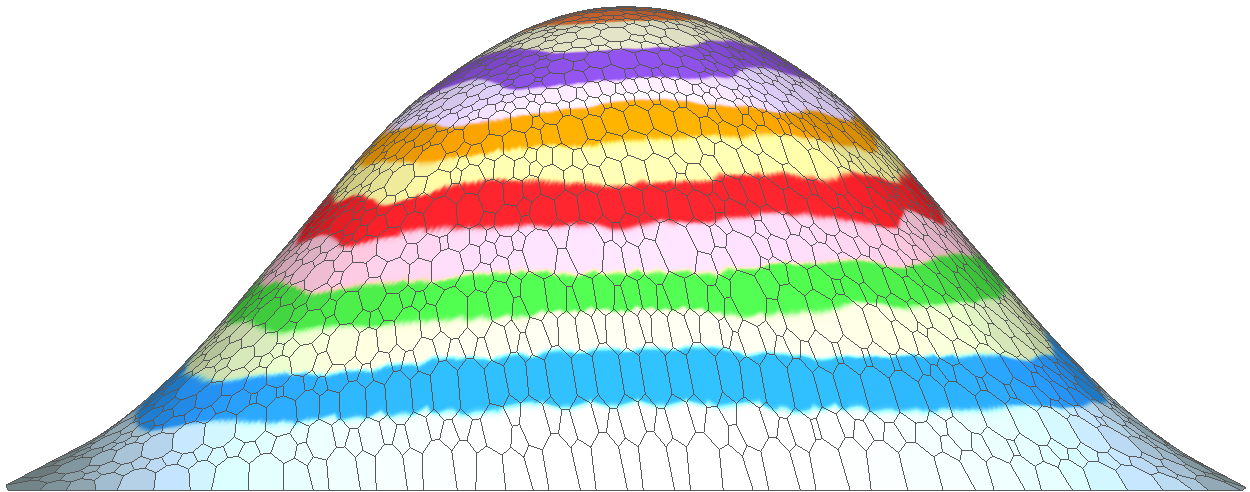}
\includegraphics[width=0.24\linewidth]{ 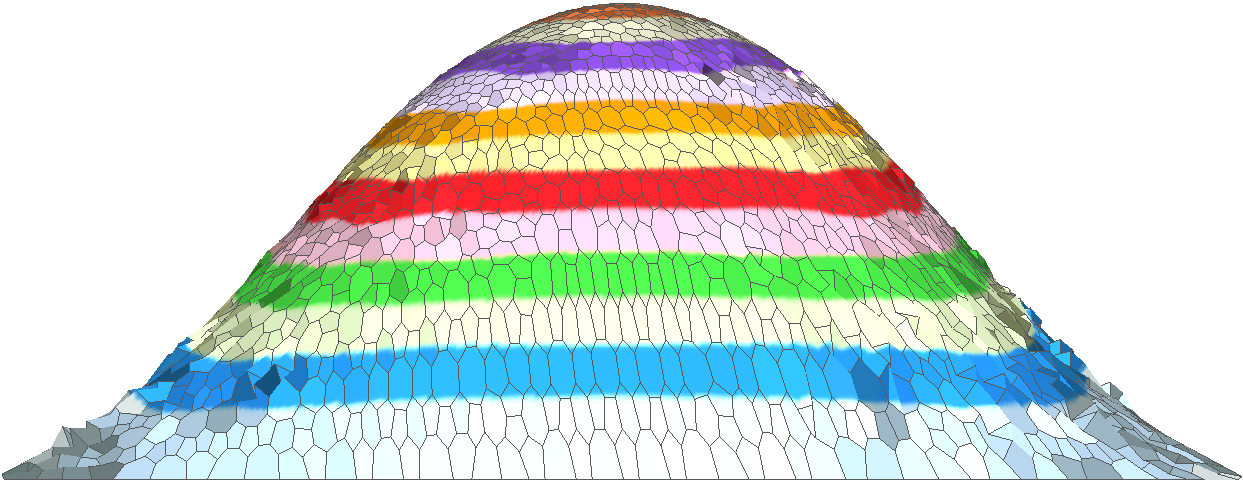}
\includegraphics[width=0.24\linewidth]{ 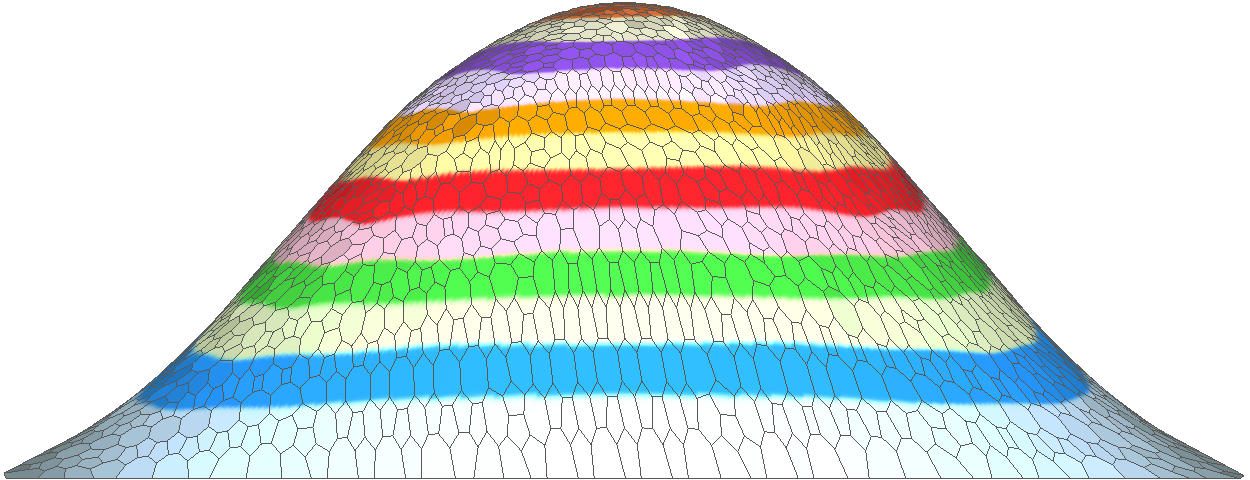}
\caption{Smoothing out bumps on a polygonal mesh (far left) with one iteration of mean curvature flow (\ref{eq:mcfWithBoundary}) with $dt=0.01$, using the Laplacians of \cite{Fujiwara1995} (center left), \cite{AlexaWardetzky2011} (center right), and \cite{PtackovaVelho2021} (farright). On general polygonal meshes, the Laplacian of \cite{PtackovaVelho2021} gives better results, as the one of \cite{AlexaWardetzky2011} quickly develops artifacts, and the Laplacian of \cite{Fujiwara1995} inflicts tangential shifts resulting in texture deformation.}\label{fig:polyMCF}
\end{figure}

In the following, we give an example of MCF of a general polygonal mesh with texture.
\begin{example}
In Figure \ref{fig:polyMCF} we illustrate the effect of MCF of a mesh formed by general polygons.
The texture mapped on the surface is a heightmap of the original mesh.
The color map highlights the fact, that $\Delta_F$ introduces a tangential component in $\vec{H}$ and, as a result, the texture drifts in the process. While the flow with $\Delta_A$ suffers from undesirable artifacts, $\Delta_P$ performs well. On the other hand, on quad meshes, engaging the Laplacians of \cite{Fujiwara1995} and \cite{AlexaWardetzky2011} results in better smoothing than with the Laplacian of \cite{PtackovaVelho2021}, in general. This is also one of the reasons we test our adaptation of DDM with Robin transmission conditions with these three different discretizations of the Laplacian --- as they all have some advantages and disadvantages, and depending on the intent, we might choose one over the other.
\end{example}

%%%%%%%%%%%%%%%%%%%%%%%%%%%%%%%%%%%%%%%%%%%%%%%%%%%%%%%%%%%%%%%%%%%%%%%%%%%

\section{Mean curvature flow with domain decomposition}\label{sec:mcfDDM}
As stated already in the introduction, we study the case of smoothing a mesh $M$ that is first divided into two connected sub-meshes $M_A$ and $M_B$, the MCF is applied on both sub-meshes with adapted Robin transmission conditions on the interface, denoted by $\Gamma$. We expect that both sub-meshes maintain a topological structure of a 2--dimensional pseudo--manifold with boundary.

\begin{figure}[ht!]
\centering
\includegraphics[width=0.3\textwidth]{ 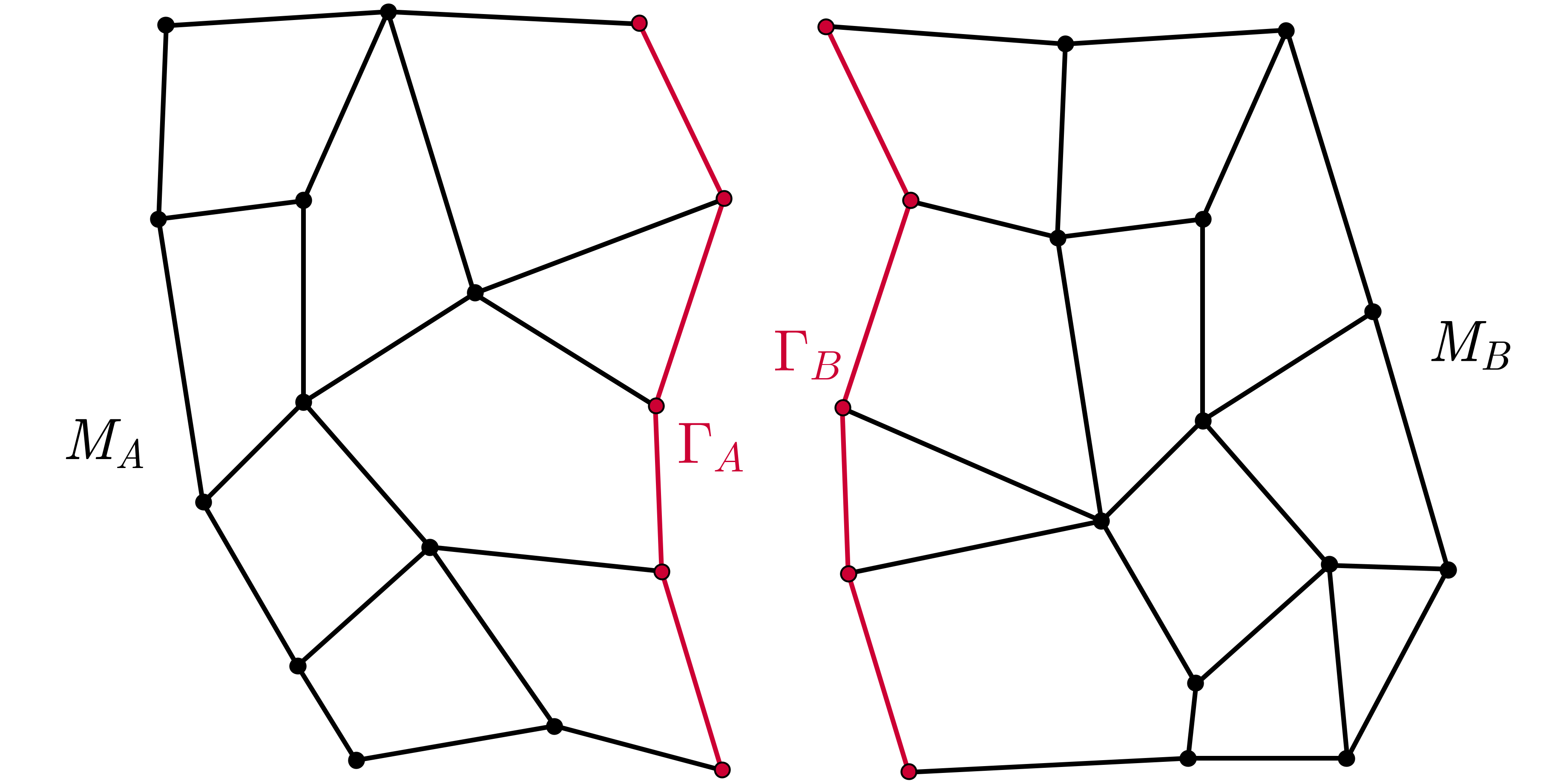}
\caption{The decomposition of a mesh $M$ into two sub-meshes $M_A$, $M_B$ without overlap. We denote the interface as $\Gamma$. The points on the interface are stored in two copies: $\Gamma_A$ belonging to $M_A$  and $\Gamma_B$ belonging to $M_B$.
}\label{fig:decomposition}
\end{figure}

\begin{figure}[ht!]
\centering
\includegraphics[width=0.4\textwidth]{ 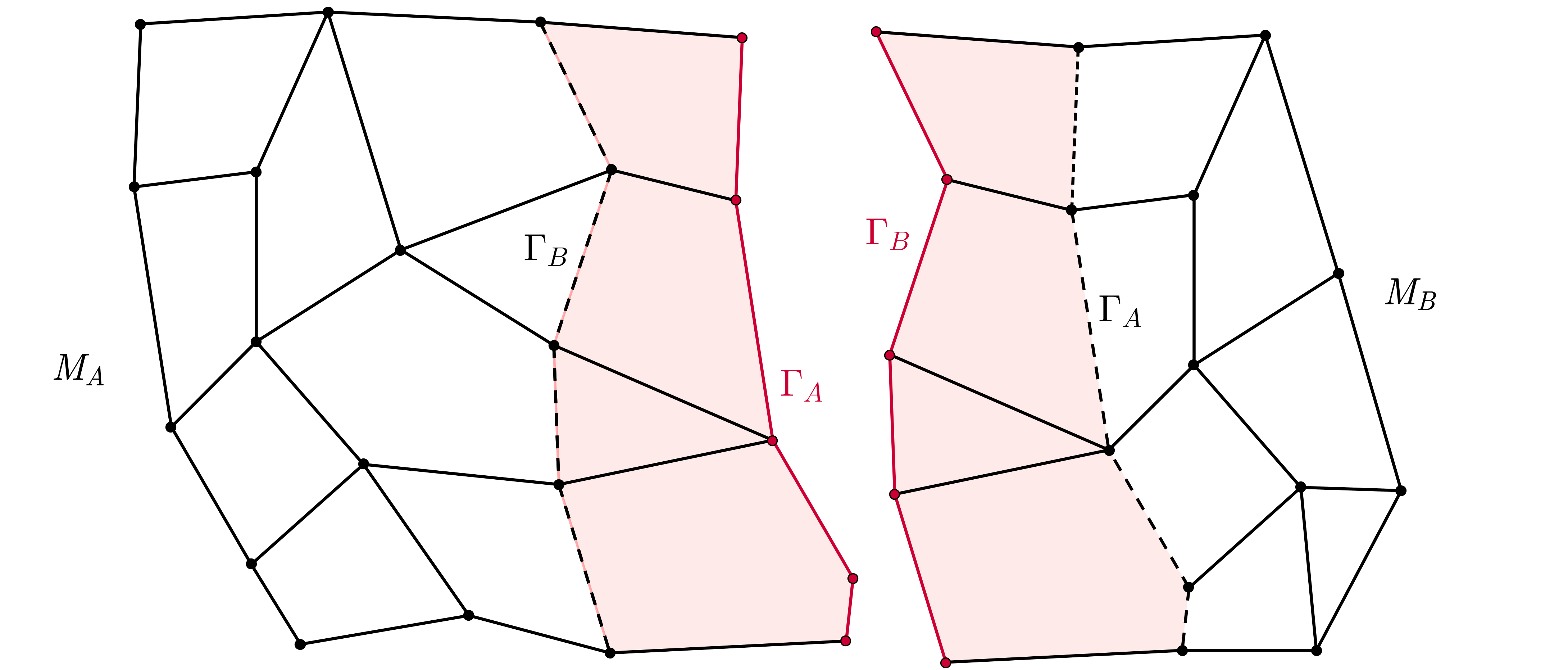}
\caption{The decomposition of a mesh into two sub-meshes $M_A$, $M_B$ with an one-face-wide overlap. The new boundary curve $\Gamma_A$ of $M_A$ corresponds now to an interior curve in $M_B$ (dashed line), and the points on boundary $\Gamma_B$ of $M_B$ correspond to a set of interior points of $M_A$.
}\label{fig:decomposition_overlap}
\end{figure}

In Figure \ref{fig:decomposition} we illustrate the decomposition of a mesh without overlap. And in Figure \ref{fig:decomposition_overlap} is shown a decomposition with an one-face-wide overlap. The decomposition without overlap is used together with Robin and Ventcell transmission conditions, while the decomposition with overlap is used in Schwarz alternating method.

We test three types of domain decomposition methods with various transmission conditions:
\begin{enumerate}
 \item Schwarz alternating method with a subdomain overlap and Dirichlet boundary conditions on the interface.
 \item Optimized Schwarz method (without overlap) with Robin and Ventcell transmission conditions.
 \item Domain decomposition without overlap with adapted Robin transmission conditions --- our actual contribution.
\end{enumerate}

%------------------------------------------------------------------------
\subsection{Schwarz alternating method with overlap}\label{subsec:SchwarzAlternating}
\begin{figure}[ht!]
\centering
\includegraphics[width=0.24\linewidth]{ 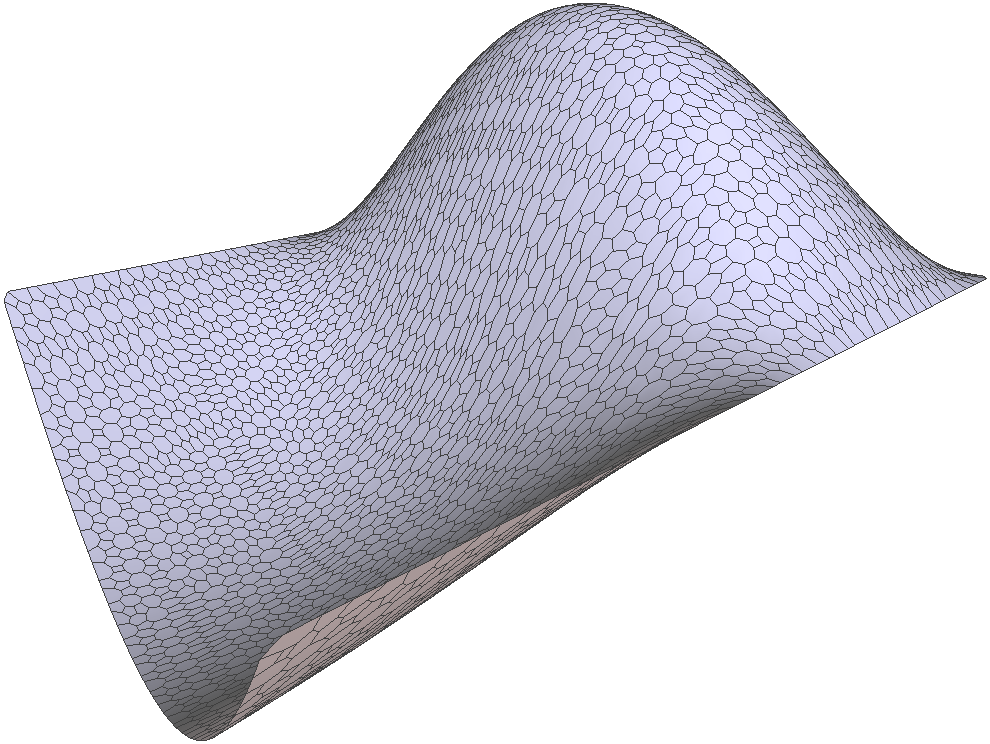}
\includegraphics[width=0.24\linewidth]{ 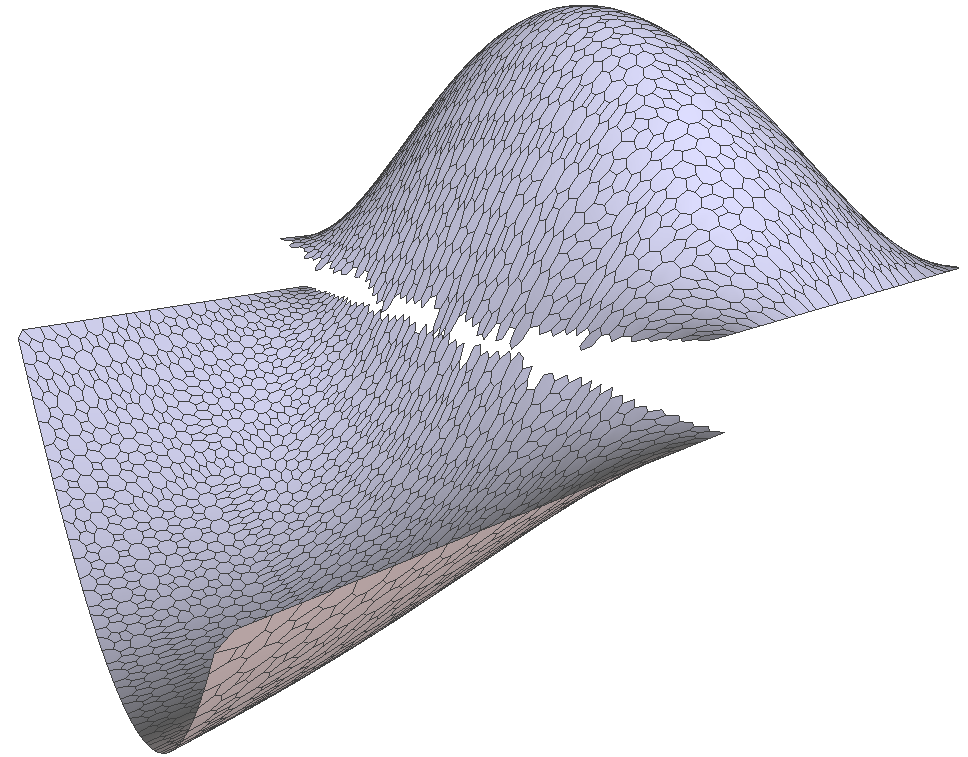}
\includegraphics[width=0.24\linewidth]{ 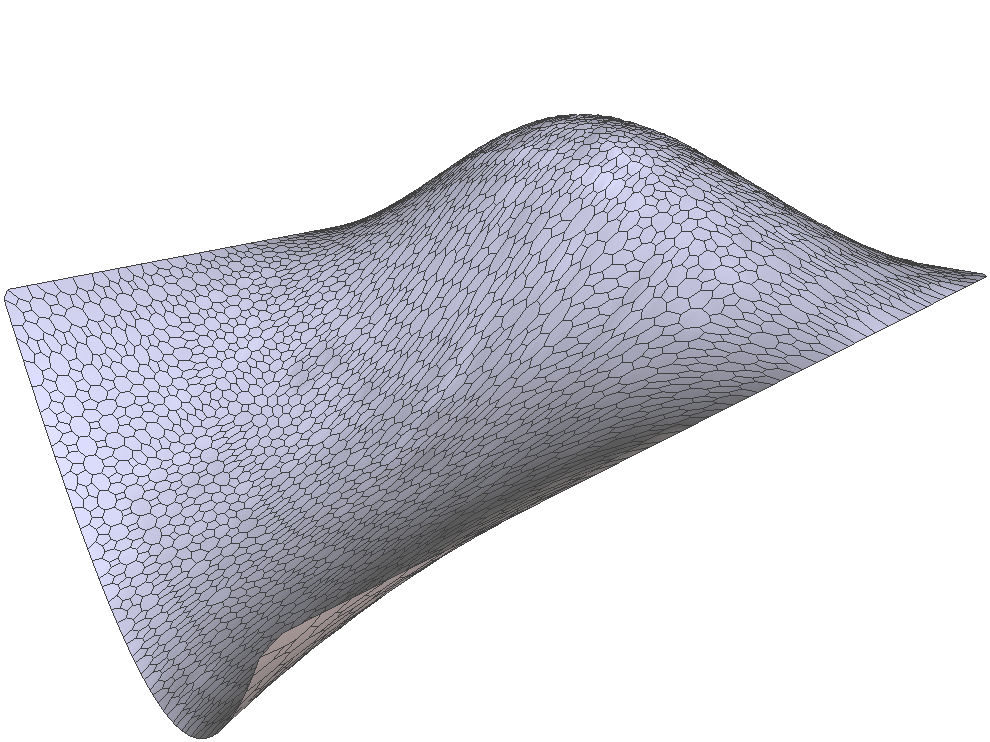}
\includegraphics[width=0.24\linewidth]{ 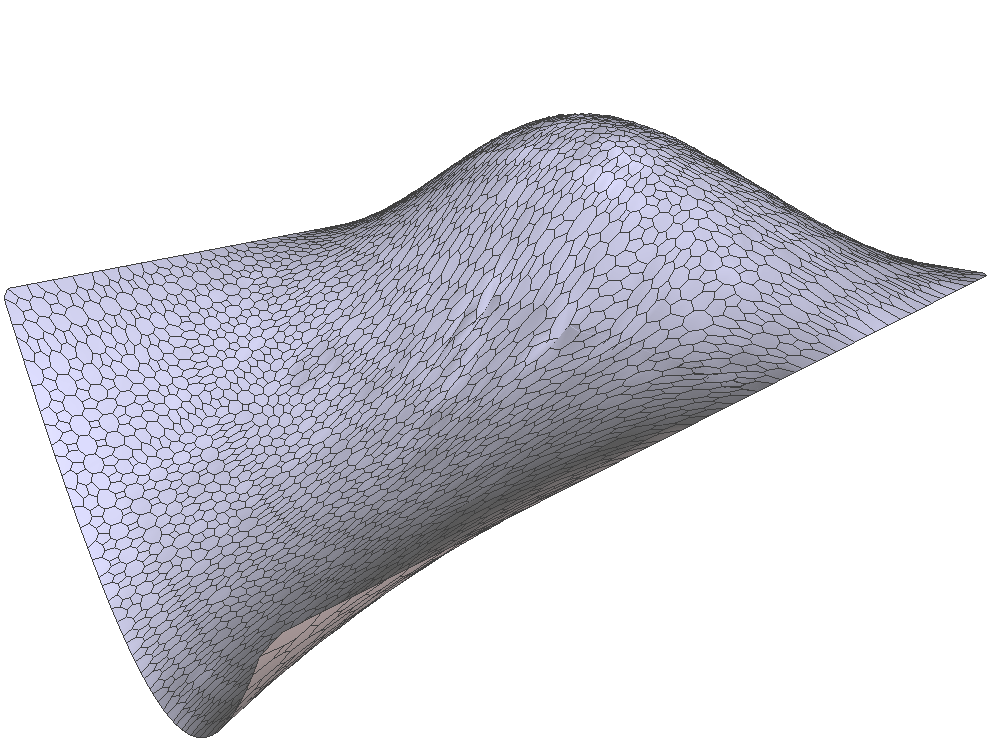}
\caption{MCF of a general polygonal mesh (upper left), which has been decomposed into two sub-meshes with an overlap (upper right). We apply two iterations of MCF with $dt = 0.05$ and employ the Laplacian of \cite{PtackovaVelho2021}. The resulting mesh without decomposition is shown in the lower left and with decomposition and Schwarz alternating method in the lower right.
}\label{fig:decomposition_poly}
\end{figure}

It is well known, that for solving the Poisson equation, the alternating Schwarz method converges slower than optimized Schwarz methods with Robin or Ventcell transmission conditions, see \cite{Gander2006}. However, there are cases, where this method comes in handy, for example if the mesh is a very irregular general polygonal mesh and the interface is not placed on a peak. See the example in the Figure \ref{fig:decomposition_poly}. We give the algorithm bellow:

Denote $V_A$ the vertices of a sub-mesh $M_A$ and $V_B$ the vertices of $M_B$, denote $\Gamma_A$ the interface vertices of $V_A$ that correspond to some interior vertices of $M_B$, see Figure \ref{fig:decomposition_overlap}, and similarly for
$\Gamma_B$. The \textbf{Schwarz alternating method with overlap} is given by these equations with boundary conditions:
\begin{equation}\label{eq:Schwarz_A}
\begin{cases}
\big(I + dt\Delta) V_A^{k+1} = V_A^k ,  & V_A\in M_A\setminus \partial M_A, \\
 V_A^{k+1} = V^{k}_A,   & V_A\in\partial M_A\setminus\Gamma_A,\\
 V_A^{k+1} = V^{k}_B,   & V_A\in\Gamma_A, V_B\in M_B\cap\Gamma_A,
\end{cases}
\end{equation}
where the vertices from $\Gamma_A\subset\partial M_A$ were set to be equal to the corresponding interior vertices of $M_B$, see the Figure \ref{fig:decomposition_overlap}, and
\begin{equation}\label{eq:Schwarz_B}
\begin{cases}
\big(I + dt\Delta) V_B^{k+1} = V_B^k ,  & V_B\setminus \partial M_B, \\
V_B^{k+1} = V^{k}_B,   & V_B\in\partial M_B\setminus\Gamma_B,\\
V_B^{k+1} = V^{k+1}_A,  & V_B\in\Gamma_B, V_A\in M_A\cap\Gamma_B,
\end{cases}
\end{equation}
where the vertices from $\Gamma_B\subset\partial M_B$ were set to be equal to the corresponding interior vertices of $M_A$, computed as solutions of (\ref{eq:Schwarz_A}).
We call the iteration complete if we solve the set of equations of (\ref{eq:Schwarz_A}) and subsequently solve the linear equations of (\ref{eq:Schwarz_B}). We thus alternate between flowing $M_A$ and $M_B$.

In Figure \ref{fig:decomposition_poly} we illustrate the MCF of a general polygonal mesh with and without decomposition. In this case, it is not possible to apply the Robin or Ventcell transmission conditions, because it is not clear how to define well the normal and tangential derivatives along the interface. However, since the interface is not placed on a peak, its vertices are not expected to move much during the MCF. Therefore, the alternating \hl{Schwarz} method performs well.

However, when we apply the standard Schwarz alternating method on a mesh and the interface is supposed to move significantly, we obtain undesirable artifacts, as illustrated in Figure \ref{fig:quad_cilinder_wavy_10200}. There, the vertices on the interface are moving much slower to the (eventually) minimal surface than the interior vertices and the mesh does not flow proportionally.

%------------------------------------------------------------------------
\subsection{Optimized Schwarz method with Robin and Ventcell transmission conditions}\label{subsec:OptimizedSchwarz}
Ventcell transmission conditions make use of normal and tangential derivatives, while Robin conditions employ only the normal derivative. With a slight abuse of proper order, we first state the conditions and then in Section \ref{subsec:NormalDerivatives} we define the discrete normal derivative and tangential derivatives.

Let $\partial_n$ denote the discrete normal derivative (Definition \ref{def:newNormalDerivative}) and $\partial_t$ the discrete tangential derivative (\ref{eq:tangentialDerivative}). Let further denote $V_A$ vertices of $M_A$ and $V_B$ vertices of $M_B$, we expect the sub-meshes to be incident only on the initial interface $\Gamma$, as in Figure \ref{fig:decomposition}. If we were to blindly rewrite the classical \textbf{Ventcell transmission conditions}\hl{ }\cite[eq. (1)]{Gander2022}\hl{ }for the mean curvature flow, they would read
\small
\begin{equation}\label{eq:Ventcell_A}
\begin{cases}
\big(I + dt\Delta) V_A^{k+1} = V_A^k ,  &  M_A\setminus \partial M_A, \\
 V_A^{k+1} = V^{k}_A,   & V_A\in\partial M_A\setminus\Gamma_A,\\
\big(I + dt(p \partial_n + q \partial_t)\big) V_A^{k+1} = \big(I + dt(p \partial_n + q \partial_t)\big) V^{k}_B,   & V_A\in\Gamma_A, V_B\in\Gamma_B,
\end{cases}
\end{equation}
\begin{equation}\label{eq:Ventcell_B}
\begin{cases}
\big(I + dt\Delta) V_B^{k+1} = V_B^k ,  &  M_B\setminus \partial M_B, \\
 V_B^{k+1} = V^{k}_B,   & V_B\in\partial M_B\setminus\Gamma_B,\\
\big(I + dt(p \partial_n + q \partial_t)\big) V_B^{k+1} = \big(I + dt(p \partial_n + q \partial_t)\big) V^{k+1}_A,   & V_A\in\Gamma_A, V_B\in\Gamma_B.
\end{cases}
\end{equation}
\normalsize
We call the iteration complete if we solve again both the equations (\ref{eq:Ventcell_A}) and subsequently (\ref{eq:Ventcell_B}).
Note that in equation (\ref{eq:Ventcell_B}) we use vertices of $\Gamma_A$ which have been computed as a solution of (\ref{eq:Ventcell_A}).
The classical \textbf{Robin transmission conditions} are the conditions of equations (\ref{eq:Ventcell_A} -- \ref{eq:Ventcell_B}) for $q= 0$.

The Ventcell (and Robin) transmission conditions with appropriate parameters $p,q$ are guaranteed to converge to the solution of the Poisson equation, and therefore find the minimal surface if the mesh is grid--like and the system is formulated as in \cite[eq. (1)]{Gander2022}. For optimal parameters and proof of convergence in such specific cases, see \cite{Gander2006, bennequin2016}. However, the conditions do not guarantee the convergence of MCF. Moreover, they do not even guarantee that the sub-meshes will be incident along the interface $\Gamma$ after each iteration, for any non--zero parameter $p$ and/or $q$. As a result, $M_A\cup M_B$ is not connected in general. Therefore, an adaptation of (\ref{eq:Ventcell_A} -- \ref{eq:Ventcell_B}) is necessary and we provide it in Section \ref{subsec:OurAlgorithm}. But we first provide the definitions of tangential and normal derivatives, as promised.

\subsubsection{Discrete normal and tangential derivatives}
\label{subsec:NormalDerivatives}
\begin{figure}[ht!]
\centering
\includegraphics[width=0.35\textwidth]{ 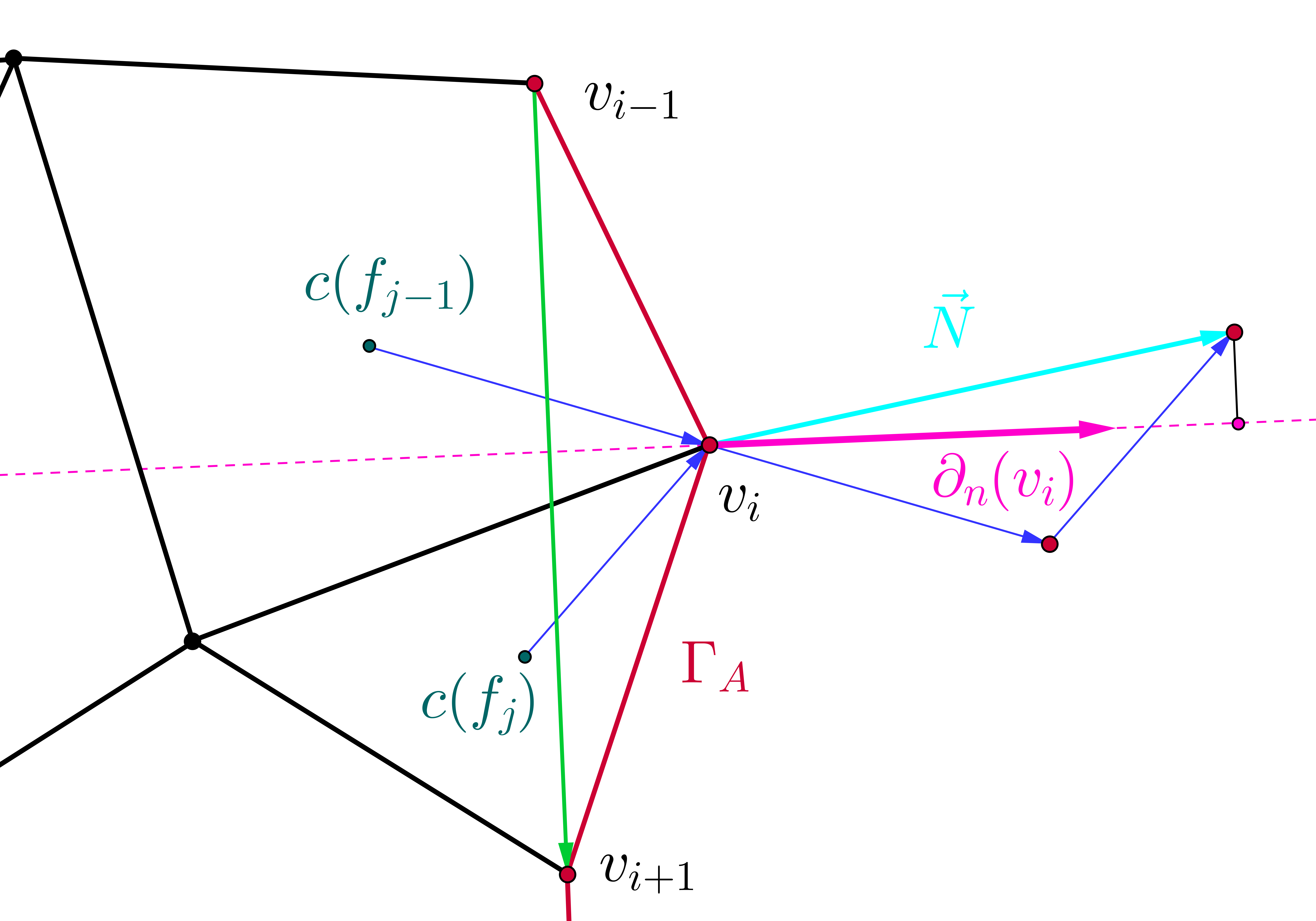}
\caption{The discrete normal derivative $\partial_n$ to the interface $\Gamma_A$ (red) at a vertex $v_i$ is a vector $\partial_n(v_i)$ colored magenta. It corresponds to the normalized orthogonal projection of vector $\vec{N}$ onto a plane, whose normal vector is $(v_{i+1} - v_{i-1})$ (green). The vector $\vec{N}$ is the sum of vectors $(v_i - c(f_{j-1}))$, $(v_i - c(f_{j}))$, where points $c(f_{j-1})$, $c(f_{j})$ are the centroids of faces incident to the vertex $v_i$.
}\label{fig:normalDerivative}
\end{figure}

The discrete normal derivative is a discrete version of directional derivative on a mesh $M$ in the direction that is outwardly 'normal' to the boundary $\Gamma$. As we can see in Figure \ref{fig:normalDerivative}, the direction perpendicular to a boundary curve at a vertex $v_i$ is not even well defined in the sense of a continuous directional derivative. We thus present its approximation in the following definition.

\begin{definition}\label{def:newNormalDerivative}
Let $v_{i-1}$, $v_{i}$, $v_{i+1}$ be three consecutive vertices on the boundary $\Gamma$ of a given mesh $M$. Denote $c(f_j)$ the centroids of faces $f_j$ incident to $v_i$, i.e., $c(f_j) = \frac{1}{p_j}\sum_{v_k\prec f_j}v_k$ for a $p_j$--gon $f_f$.
Further, denote
$\vec{t} = \frac{v_{i+1}-v_{i-1}}{\norm{v_{i+1}-v_{i-1}}_2}$.
The \textbf{discrete normal derivative} $\partial_n$ to $\Gamma$ at vertex $v_i$ is:
\begin{equation}
\begin{array}{rcl}
\partial_n(v_i) &=&  \frac{\vec{N}}{\norm{\vec{N}}_2}, \text{ where}\\
\vec{N} &=& \displaystyle\sum_{f_j\succ v_i} (v_i - c(f_j)) -
\big\langle \displaystyle\sum_{f_j\succ v_i} (v_i - c(f_j)) \;,\;
\vec{t}\;\big\rangle\;
\vec{t}.
\end{array}
\end{equation}
\end{definition}
On a planar uniform grid, the above definition would correspond to the normal derivative of standard second order five point finite difference discretization as in \cite[Section 2]{Gander2022}.

In our implementation, we can easily compute the centroids of faces using the matrix $\fv$ of eq. (\ref{eq:fv}):
\[
\begin{pmatrix}
c(f_0)\\
\vdots \\
c(f_{|F| - 1})
\end{pmatrix}=
\Big(
\diag\big((1,\dots,1) \lceil \fv \rceil \big)
- \lceil \fv^\top \rceil \fv
\Big)
\begin{pmatrix}
V_0\\
\vdots \\
V_{|V| - 1},
\end{pmatrix}
\]
where the ceiling function $\lceil \cdot \rceil$ is applied element--wise.

Let again $v_{i-1}$, $v_{i}$, $v_{i+1}$ be three consecutive vertices on the interface $\Gamma$. The \textbf{discrete tangential derivative} $\partial_t$, i.e., directional derivative along the interface $\Gamma$ at vertex $v_i$ is given by:
\begin{equation}\label{eq:tangentialDerivative}
\partial_t(v_i) =
\frac{2}{\norm{v_{i} - v_{i-1}} + \norm{v_{i} - v_{i+1}}}
\Big( \frac{v_{i} - v_{i-1}}{\norm{v_{i} - v_{i-1}}} +
\frac{v_{i} - v_{i+1}}{\norm{v_{i} - v_{i+1}}}  \Big),
\end{equation}
where the norms are standard Euclidean. Once again, on an uniform grid, the tangential derivative would correspond to the standard finite difference discretization \hl{given} in \cite[Section 2]{Gander2022}.

%------------------------------------------------------------------------

\subsection{Adapted Robin and Ventcell transmission conditions}\label{subsec:OurAlgorithm}
We now adapt the Ventcell transmissiion conditions so that after each complete iteration, i.e., after flowing both parts $M_A$ and $M_B$ of the initial mesh, the reassembled mesh $M = M_A \cup M_B \setminus \Gamma_B$ is connected (and as smooth as possible).

Let $\partial_n$ denote the discrete normal derivative (Definition \ref{def:newNormalDerivative}) and $\partial_t$ the discrete tangential derivative (\ref{eq:tangentialDerivative}). Let further denote $V_A$ the vertices of $M_A$ and $V_B$ the vertices of $M_B$. Our \textbf{adapted Ventcell transmission conditions} are
\begin{equation}\label{eq:adaptedVentcell__A}
\begin{cases}
\big(I + dt\Delta) V_A^{k+1} = V_A^k ,  & V_A\in M_A\setminus\partial M, \\
 V_A^{k+1} = V^{k}_A,   & V_A\in\partial M_A\setminus\Gamma_A,\\
\big(I + dt(p \partial_n + q \partial_t)\big) V_A^{k+1} = (I - dtp \partial_n) V^{k}_B,   & V_A\in\Gamma_A, V_B\in\Gamma_B,
\end{cases}
\end{equation}
\begin{equation}\label{eq:adaptedVentcell__B}
\begin{cases}
\big(I + dt\Delta) V_B^{k+1} = V_B^k ,  & V_B\in M_B\setminus\partial M_B, \\
V_B^{k+1} = V^{k}_B,   & V_B\in\partial M_B\setminus\Gamma_B,\\
V_B^{k+1} = V^{k+1}_A,   & V_A\in\Gamma_A, V_B\in\Gamma_B.
\end{cases}
\end{equation}
We call the iteration complete if we solve the set of equations of (\ref{eq:adaptedVentcell__A}) and subsequently solve the linear equations of (\ref{eq:adaptedVentcell__B}).
Note that in equation (\ref{eq:adaptedVentcell__B}) we use vertices of $\Gamma_A$ which have been computed as a solution of (\ref{eq:adaptedVentcell__A}).
Again, the adapted \textbf{Robin transmission conditions} are the conditions of equations (\ref{eq:adaptedVentcell__A} -- \ref{eq:adaptedVentcell__B}) for $q= 0$.

\subsubsection*{Ideal parameters}
We have performed \hl{numerous} experiments with various time steps, three different discrete Laplacians, on many various meshes, just to find out that the \hl{overall} best performance with the adapted Robin transmission conditions is for the parameter $p^\ast$ given by the closed form \cite[eq. (4.15)]{Gander2006}:
\begin{equation}\label{eq:parameterP}
 p^\ast = \big( (k^2_{\min} + \eta) (k^2_{\max} + \eta)\big){}^\frac{1}{4},
\end{equation}
where the problem parameter $\eta = 1$, $k_{\min}$ is the lowest representable frequency of the mesh, and $k_{\max}$ is the maximal frequency that can be represented by the given mesh. Based on the findings of \cite[proof of Theorem 4.2]{Gander2006}, adapted to our case, we set
\[
 k_{\max} = \frac{\pi}{\norm{e}},\quad k_{\min} = \frac{\pi}{l},
\]
where $\norm{e}$ is minimal edge length of the input mesh, and $l$ is the width of the initial sub-mesh (an approximate distance between the vertices of $\Gamma$ and the vertices of the opposite boundary). Summing up, the equation (\ref{eq:parameterP}) then reads
\begin{equation}\label{eq:ourP}
 p^\ast = \Big( \Big(\Big(\frac{\pi}{l}\Big)^2 + 1 \Big)
 \Big(\Big(\frac{\pi}{\norm{e}}\Big)^2 + 1 \Big)\Big){}^\frac{1}{4}.
\end{equation}

\begin{figure}[tbp]
\centering
\includegraphics[width=0.24\textwidth]
{ 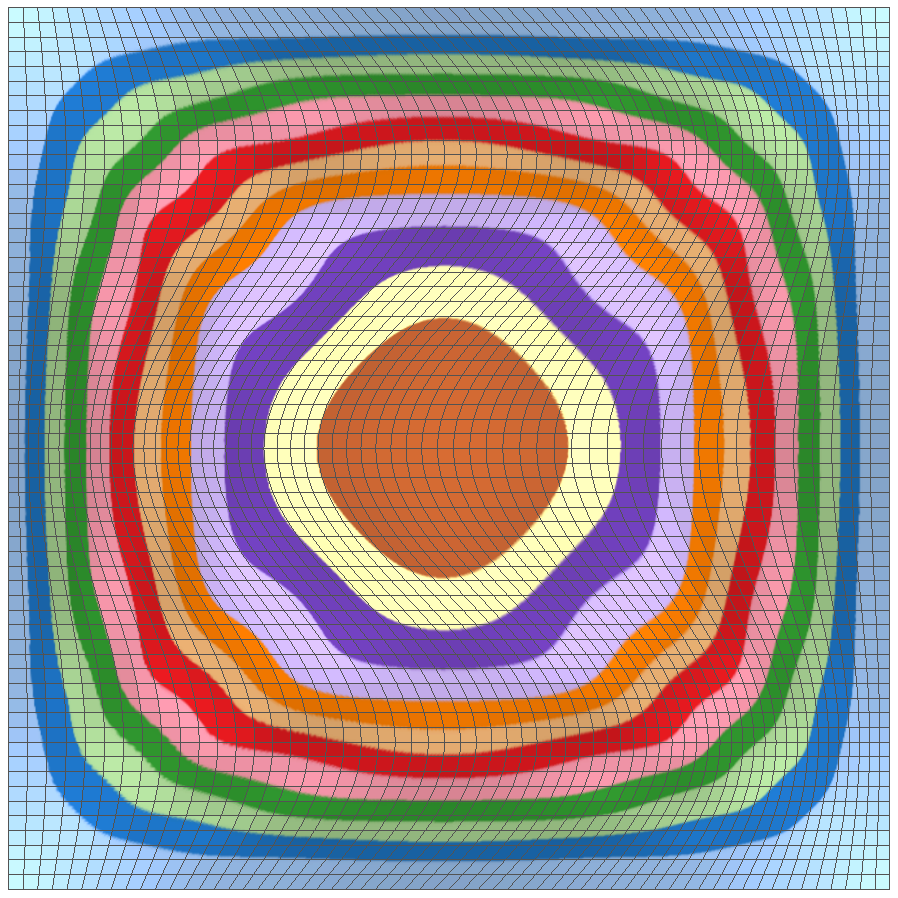}
\includegraphics[width=0.24\textwidth]
{ 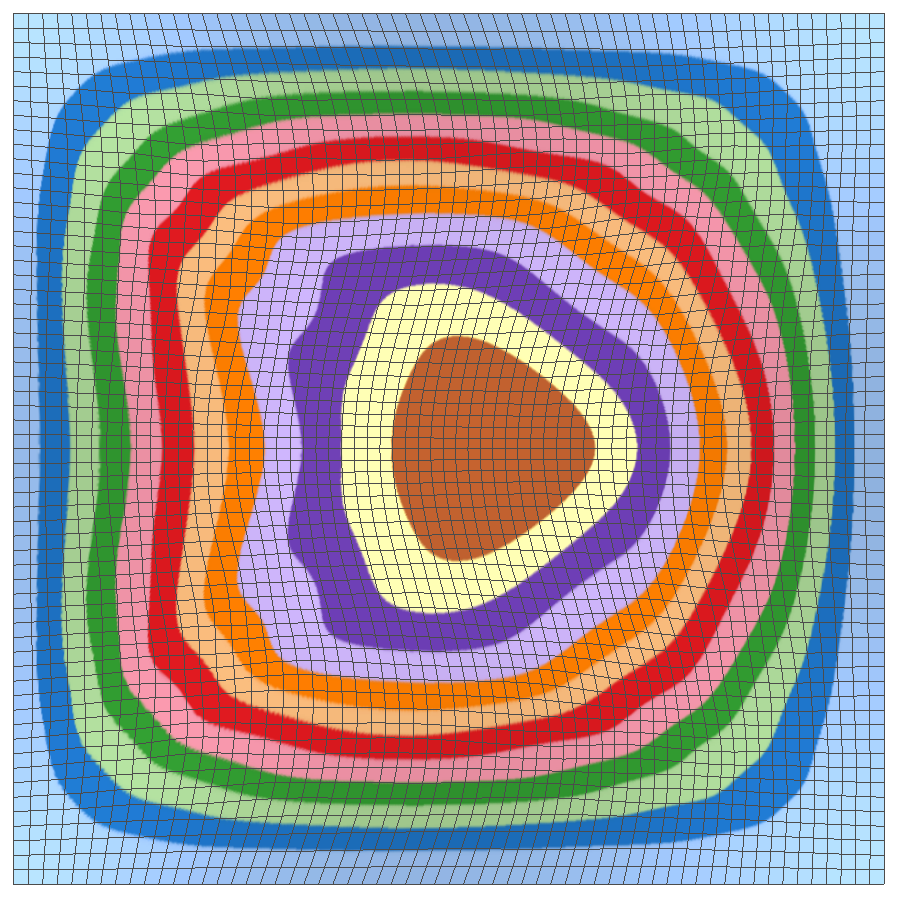}
\includegraphics[width=0.24\textwidth]
{ 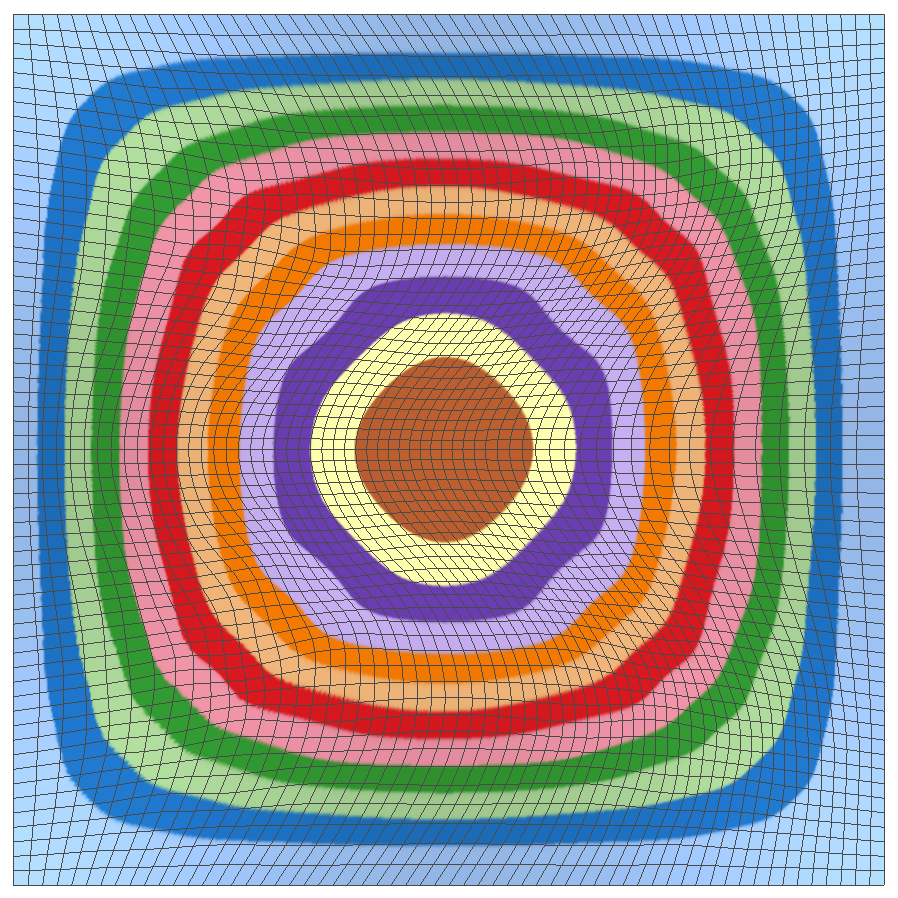}
\includegraphics[width=0.24\textwidth]
{ 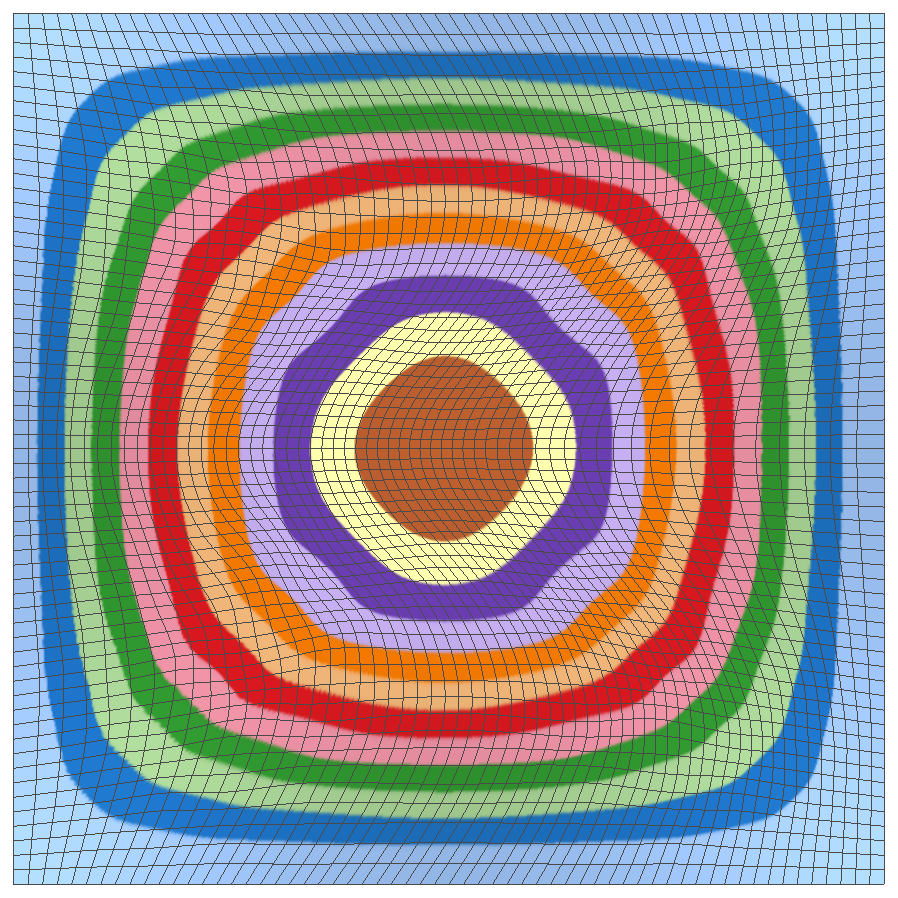}
\hfill
\\
\includegraphics[width=0.24\textwidth]
{ 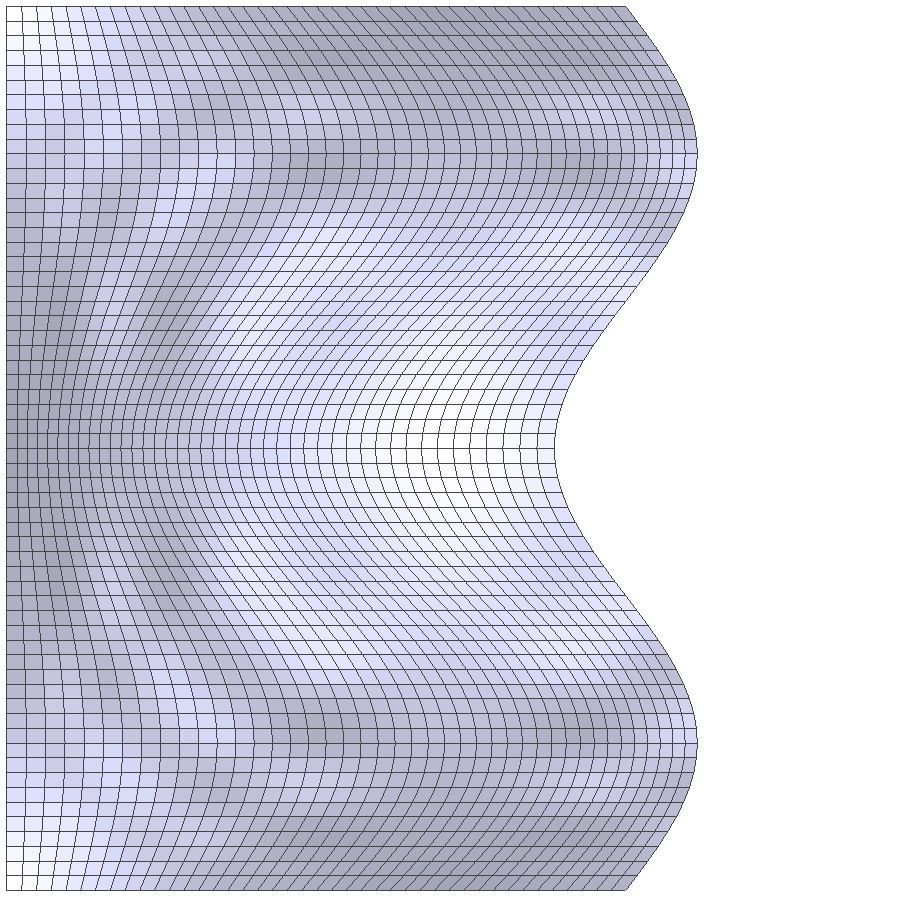}
\includegraphics[width=0.24\textwidth]
{ 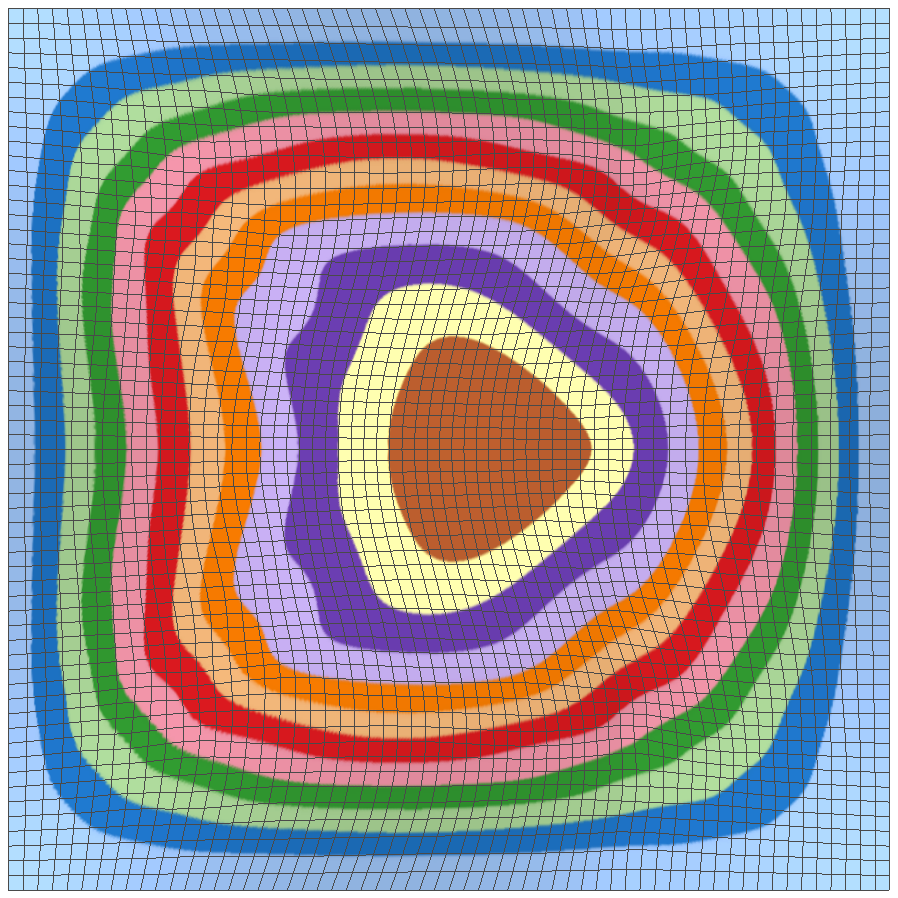}
\includegraphics[width=0.24\textwidth]
{ 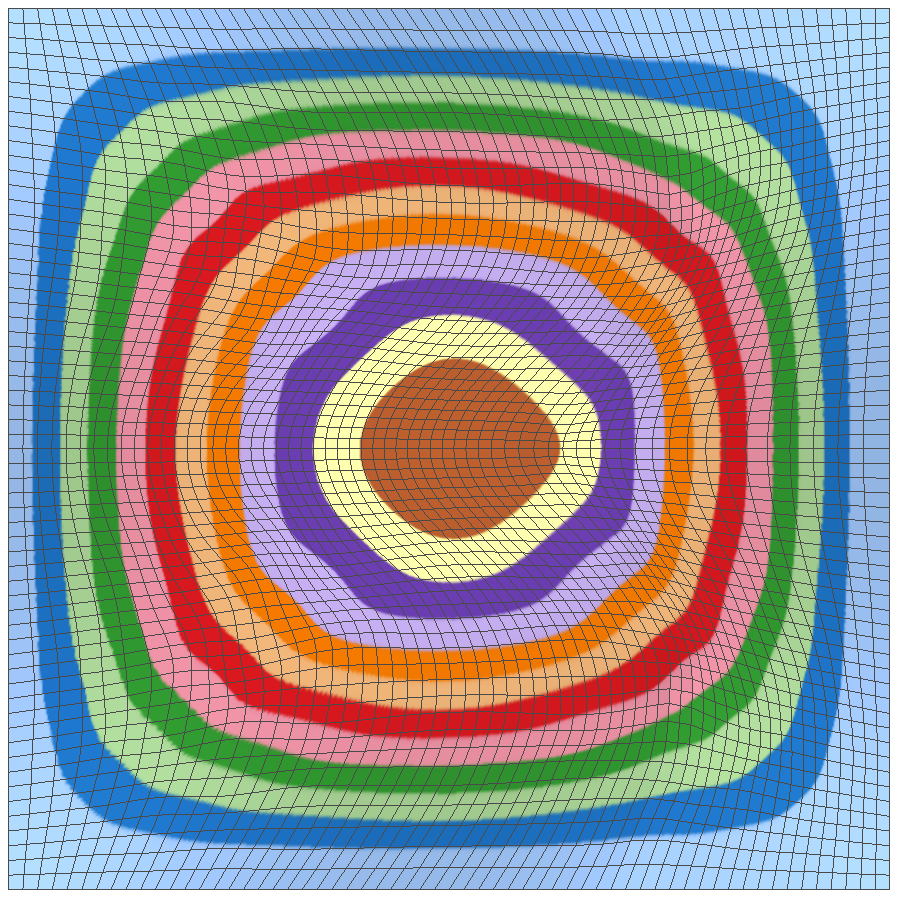}
\includegraphics[width=0.24\textwidth]
{ 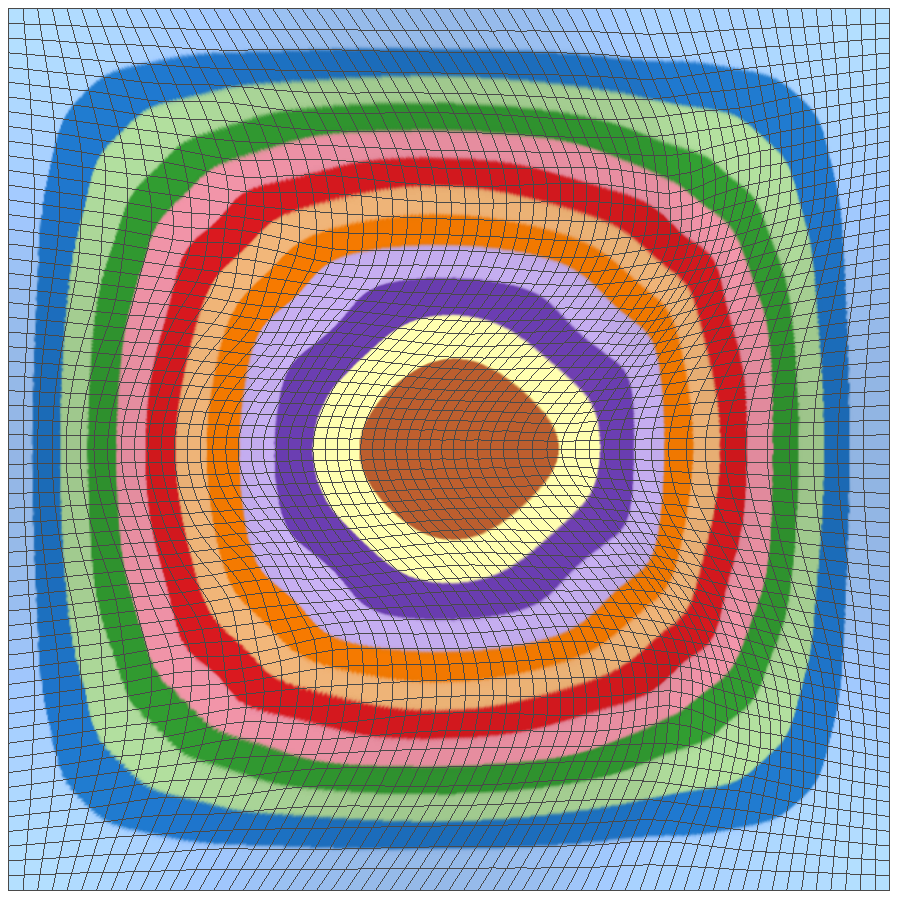}
\hfill
\caption{%
Mesh smoothing through mean curvature flow without and with domain decomposition; with time step $dt = 0.05$. In the top left corner is the original mesh with a texture of its heightmap; the underlying surface is the same as in Figure \ref{fig:polyMCF}.
Then from right to left follow meshes after 1 iteration of MCF with Laplacians of \cite{Fujiwara1995}, \cite{AlexaWardetzky2011}, and \cite{PtackovaVelho2021}.
In the bottom row we perform smoothing after first decomposing the original mesh into two sub-meshes; one of its part is the far left image. We use the Laplacians in the same order as in the upper row. When we chose appropriate conditions on the interface, it is hard to distinguish between the corresponding meshes.
}\label{fig:quad_tex_3721_top}
\end{figure}

In general, we set the parameter $q$ to zero, since the effect of employing the tangential derivative (\ref{eq:tangentialDerivative}) is the straightening of the curve on the interface. However, the interested reader is referred to \cite{Gander2006, bennequin2016} to extract the appropriate parameters $p, q$ for the Ventcell transmission conditions --- if the curve on the interface $\Gamma$ is rather straight. For example, on grid-like meshes we can employ the tangential derivative without introducing any artifacts if we perform MCF with $\Delta_F$, as in Figure \ref{fig:quad_tex_3721_top}, since the Laplacian of \cite{Fujiwara1995} itself has a straightening effect on the interface curve $\Gamma$.

\subsubsection*{Correction of vertices on the interface}
When we perform a domain decomposition and apply MCF on submeshes $M_A$, $M_B$ by solving equations (\ref{eq:adaptedVentcell__A} -- \ref{eq:adaptedVentcell__B}), we get $\Gamma_A = \Gamma_B$. However, the reassembled mesh
$M = M_A \cup M_B \setminus \Gamma_B$ can become non--smooth close to the interface $\Gamma$. We easily minimize this issue by a correction after each complete iteration.

We perform a very simple correction by substituting vertices $v_i \in\Gamma_A$ for centroids of their neighbourhood vertices, i.e., we redefine them as
\begin{equation}\label{eq:correction}
 v_i = \frac{1}{n} \sum_{v_j\sim v_i} v_j,\quad \text{where $n$ is the valence of $v_i$.}
\end{equation}

%------------------------------------------------------------------------
%------------------------------------------------------------------------

\section{Implementation and evaluation}\label{sec:Implementation}
We have implemented our algorithm for alternating Schwarz method (\ref{eq:Schwarz_A} -- \ref{eq:Schwarz_B}) and for optimized Schwarz method with adapted Ventcell transmission conditions (\ref{eq:adaptedVentcell__A} -- \ref{eq:adaptedVentcell__B}) in the Python programming language. The code will be publicly available at Github at the time of publication of this article.

All of the methods are almost fully parallelizable. We can construct the two matrices of the chosen discrete Laplacian corresponding to sub-meshes $M_A$ and $M_B$ in parallel. Constructing these matrices involves integration over mesh elements and it is therefore computationally expensive and time-consuming. For solving the set of equations (\ref{eq:Schwarz_A} -- \ref{eq:Schwarz_B}), resp. (\ref{eq:adaptedVentcell__A} -- \ref{eq:adaptedVentcell__B}), we do need to know the right--hand sides, thus we need to know the position of the vertices on the interface $\Gamma$. Therefore, solving for $M_A$ and $M_B$ can not run fully simultaneously. However, the most computation--intensive part, i.e., the LU-- or QR-- factorization of the subdomain matrices, is parallelizable.

Our implementation is not optimized for comparing the computational times. However, it is well known in the DDM community, that the speed-up of the computation by virtue of parallelizing (some parts of) is significant. Concretely, a single iteration cost amounts to the solution cost of the largest subdomain \cite[Chapter 8]{Dolean2015Intro}. Generally, this saving has to be compared with the need for multiple iterations to achieve a sufficient accuracy. This trade-off and the scalability of DDMs has been showcased in \cite[Chapter 8]{Dolean2015Intro}, for instance.
Moreover, there is a well-established link between the performance of the fully parallelizable version and the one we are using, see \cite{Gander2019}.

\begin{figure}
\centering
\includegraphics[width=0.24\linewidth]{ 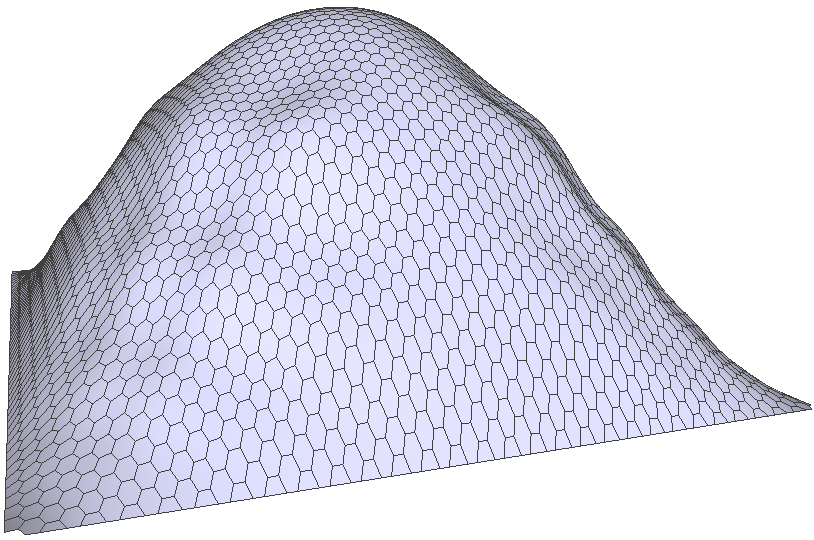}
\includegraphics[width=0.24\linewidth]{ 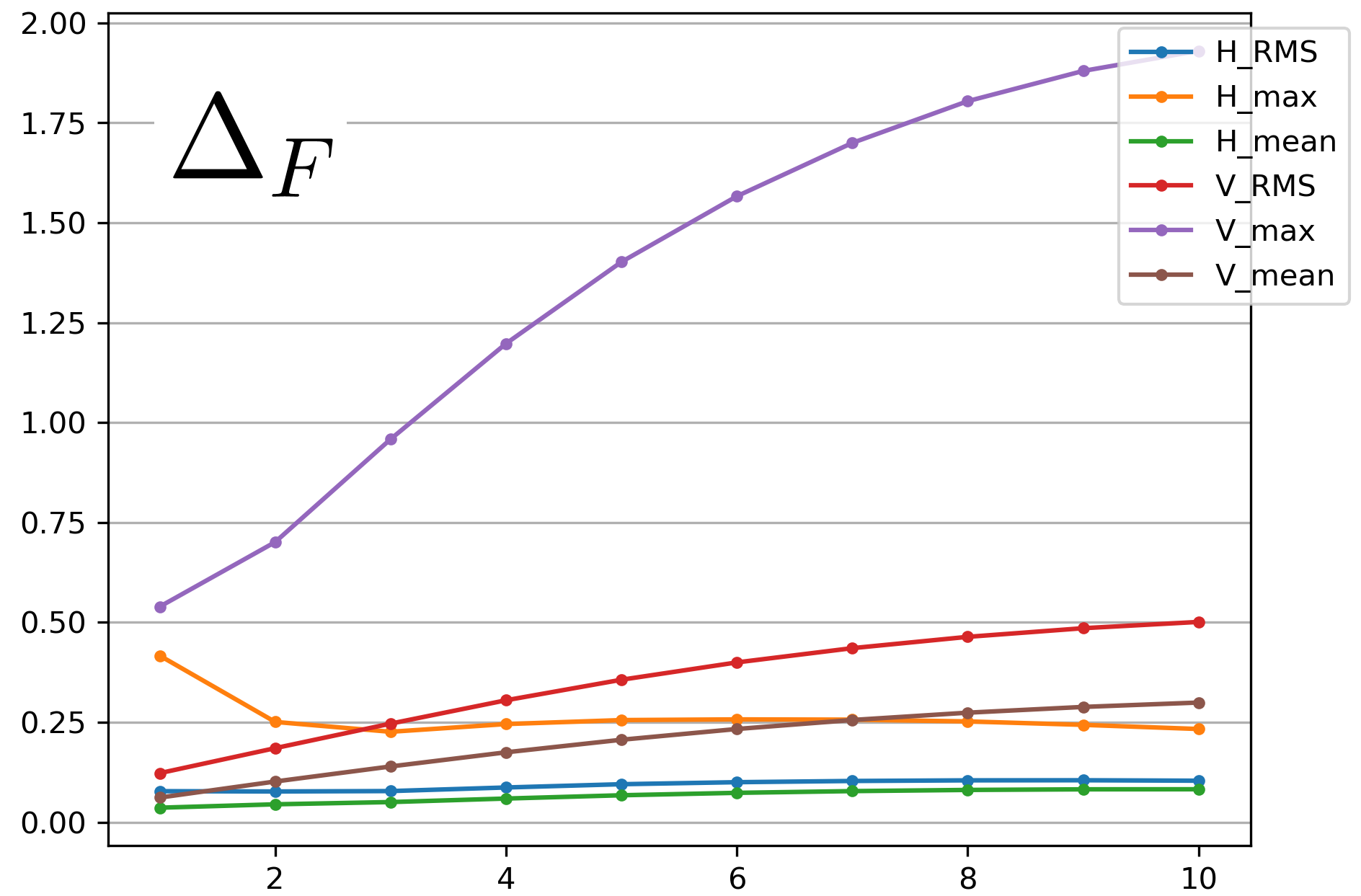}
\includegraphics[width=0.24\linewidth]{ 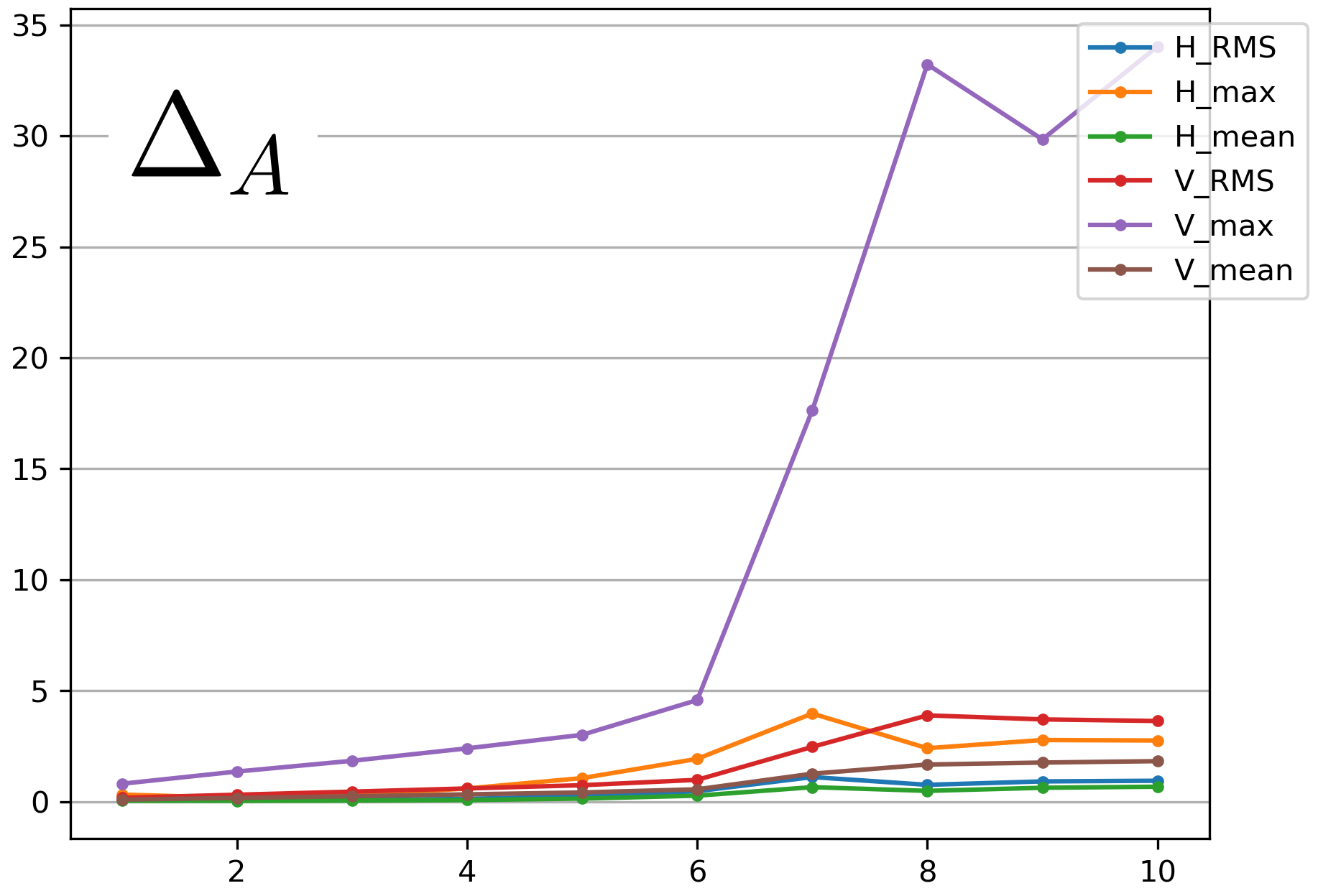}
\includegraphics[width=0.24\linewidth]{ 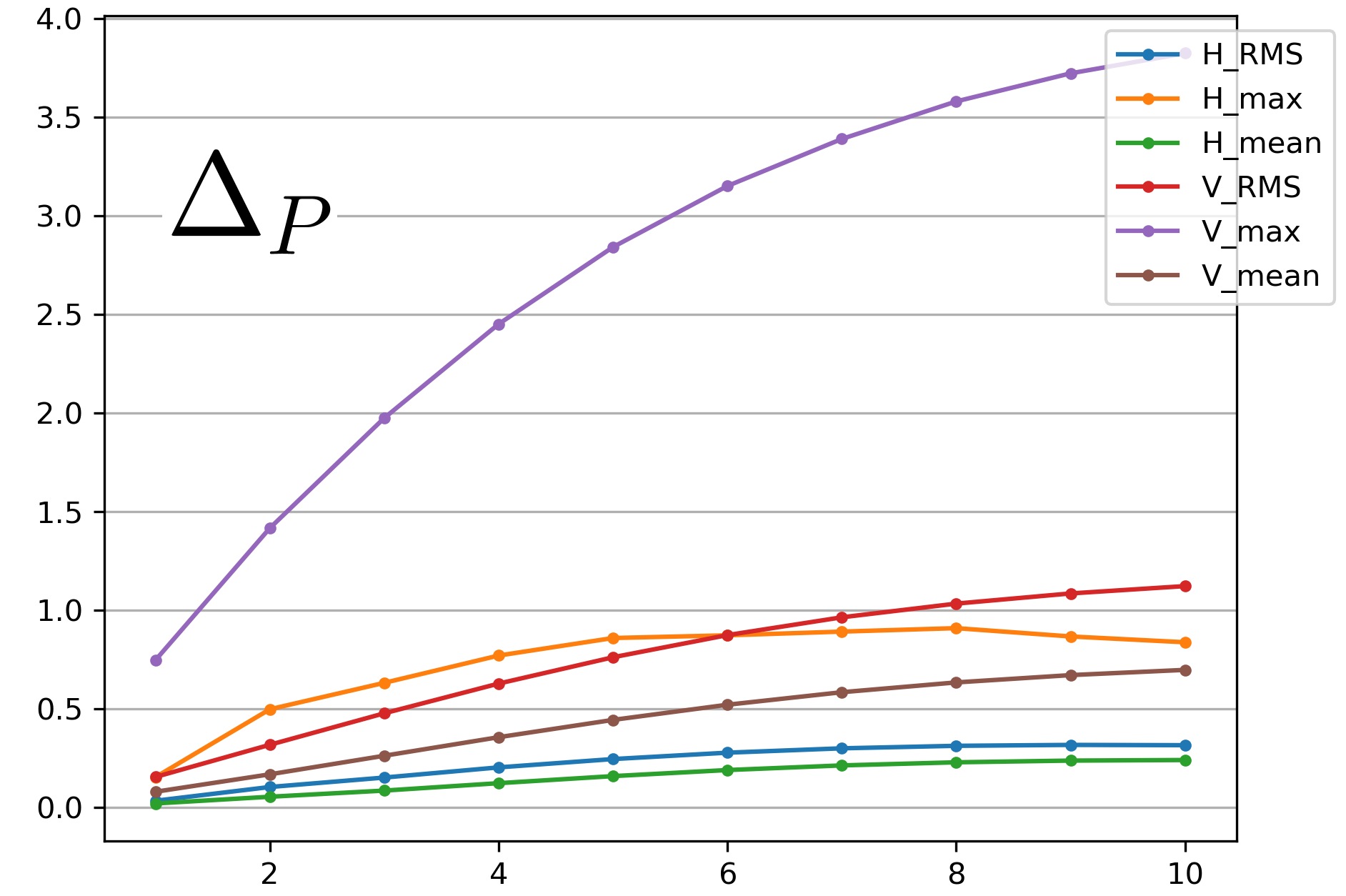}
\caption{MCF with domain decomposition on a prevalently hexagonal mesh (6118 vertices, $\min\norm{e} \approx 0.005$). In the far left is the original mesh, which is automatically divided into two approximately equal sub-meshes; $dt = 0.01$. Then follow the plots of relative distances between the smoothed reference mesh and the mesh smoothed by MCF with our adapted method (\ref{eq:adaptedVentcell__A}--\ref{eq:adaptedVentcell__B}) for $\Delta_F$, $\Delta_A$, and $\Delta_P$. The axis $y$ corresponds to the number of iterations and the axis $x$ shows the relative distances (see Section \ref{sec:Evaluation} for an explanation of the symbols in the legend). Except for the Laplacian $\Delta_A$, the distances of the corresponding meshes are small. However, for $\Delta_A$ even without decomposition the mesh flows into a non-smooth spiky mesh, starting at iteration 5.
}\label{fig:hex_F_tex_6118}
\end{figure}

\begin{figure}[h]
\centering
\includegraphics[width=0.24\textwidth]{ 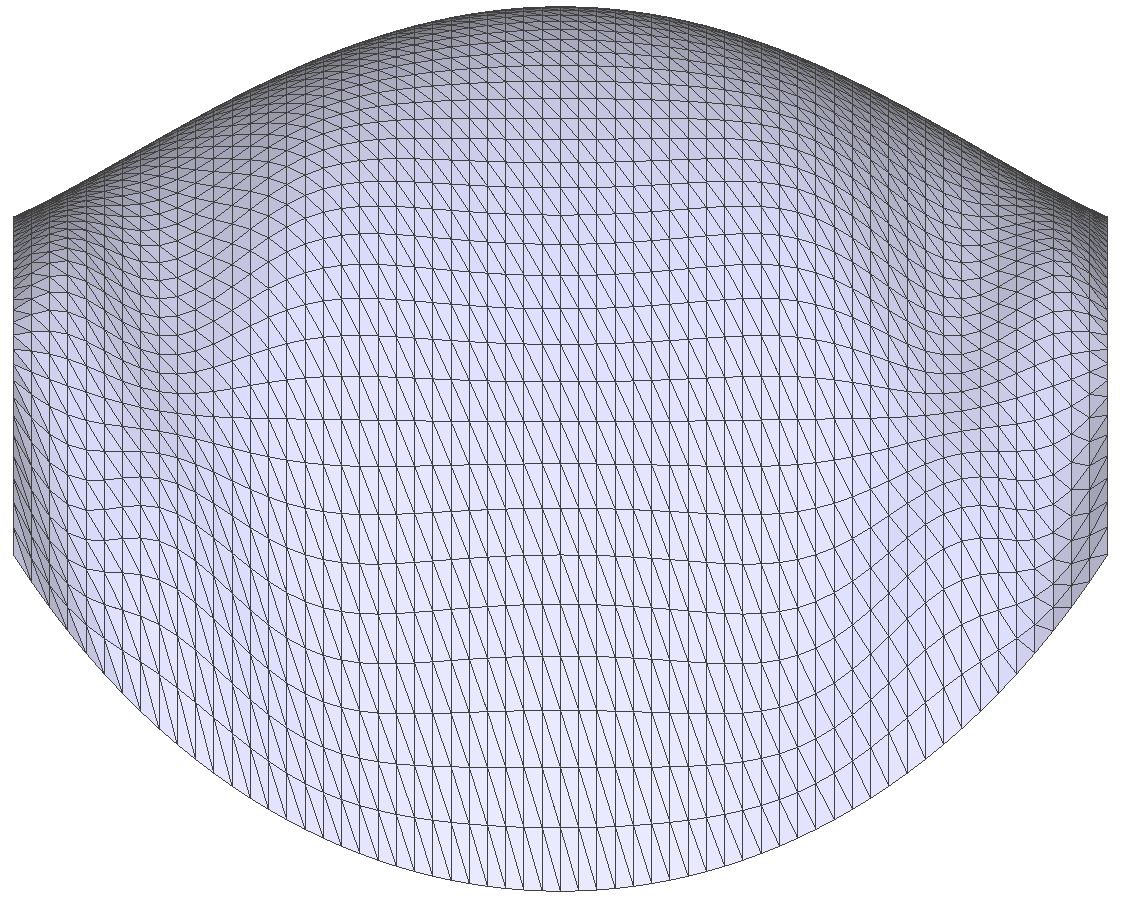}
\includegraphics[width=0.24\textwidth]{ 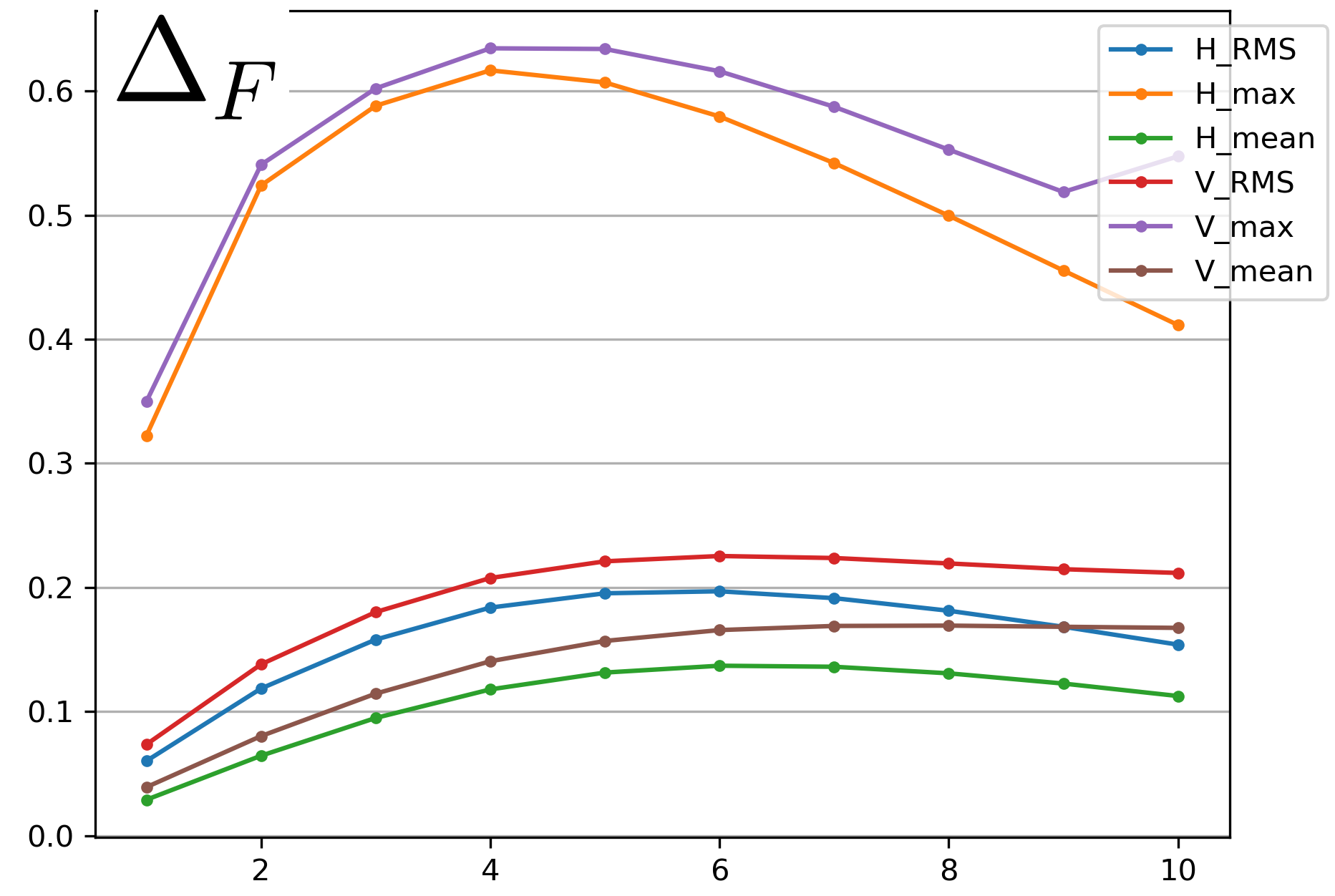}
\includegraphics[width=0.24\textwidth]{ 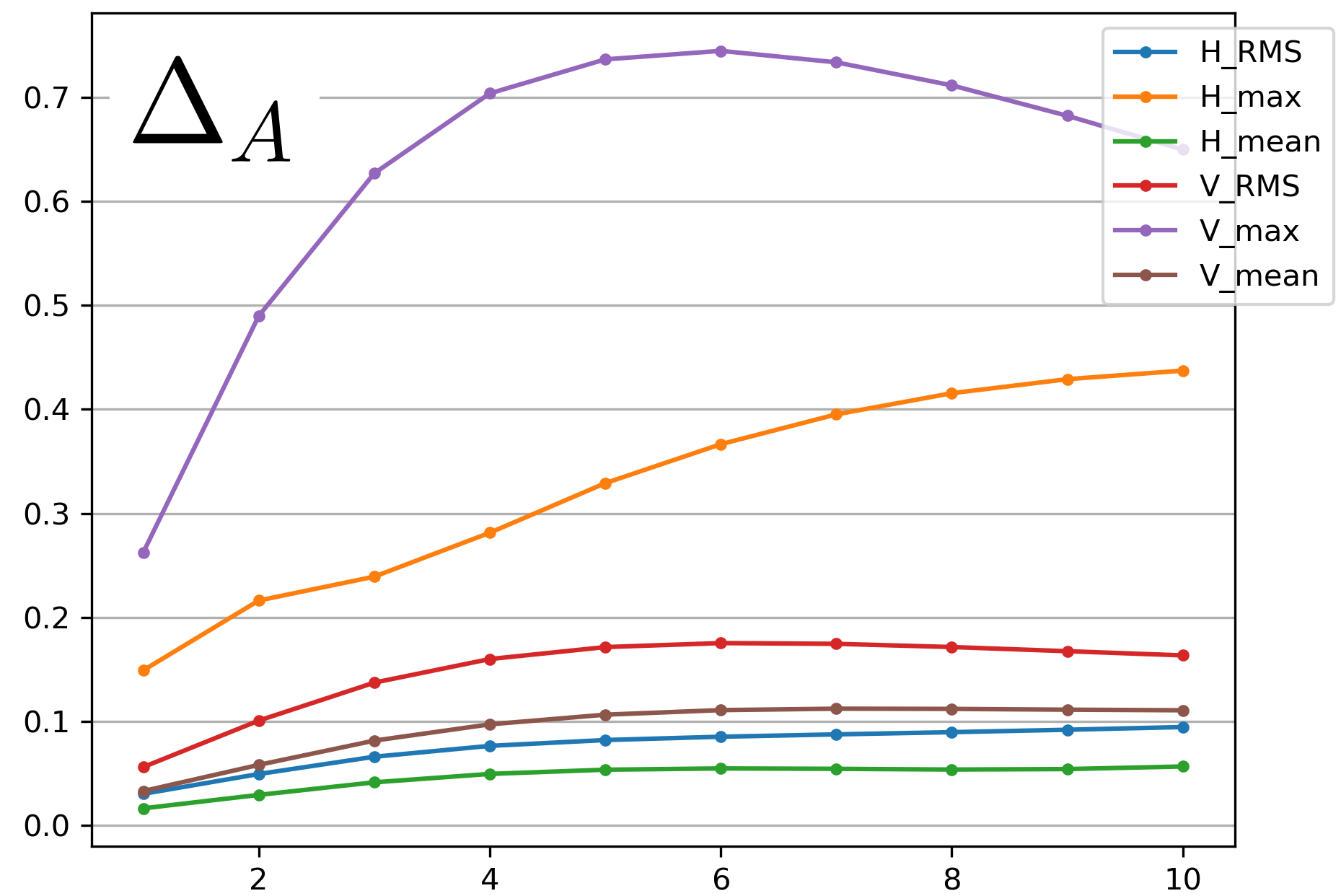}
\includegraphics[width=0.24\textwidth]{ 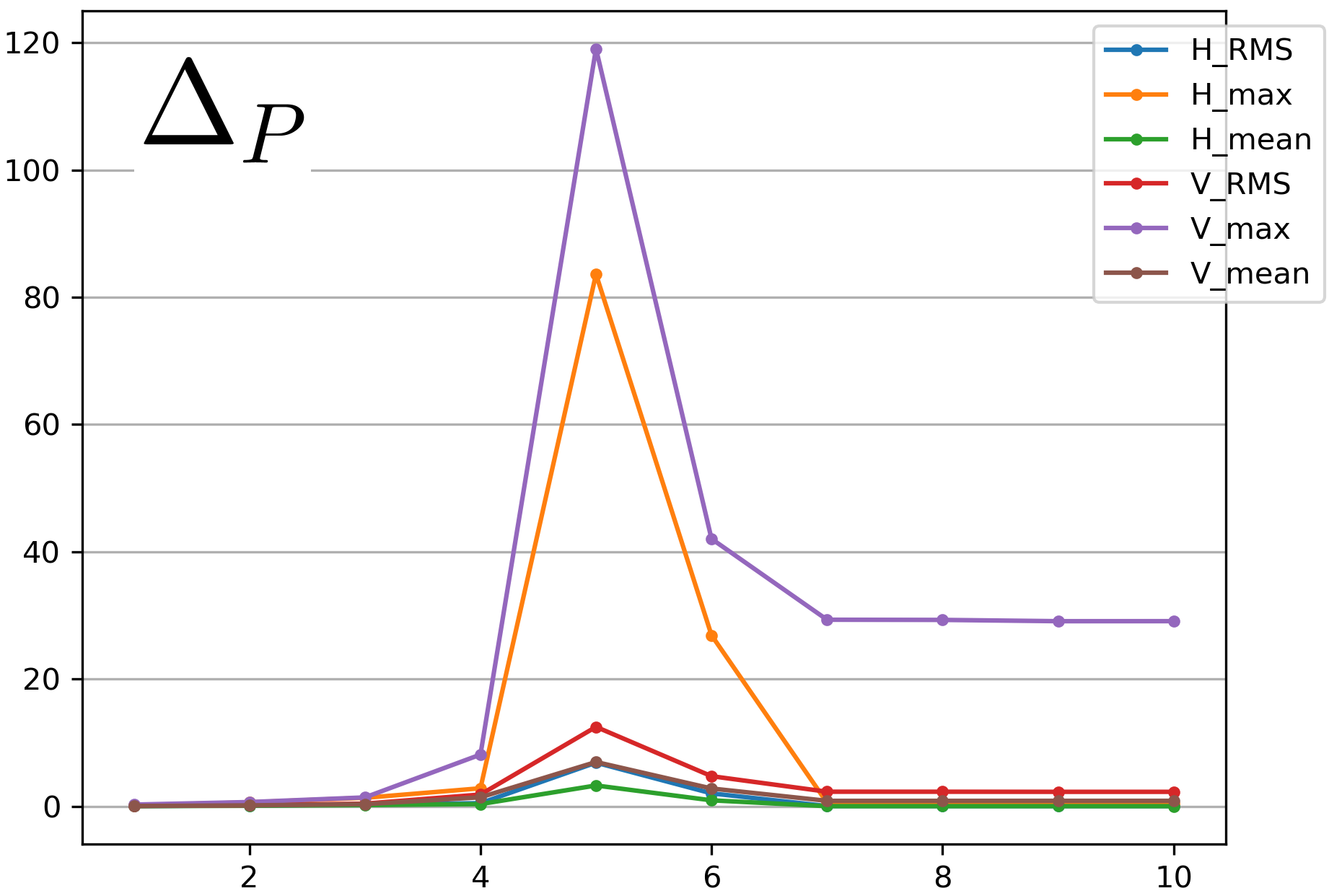}
\caption{MCF with domain decomposition of a triangle mesh (3721 vertices, $\min\norm{e} \approx 0.03$). The mesh is automatically divided into two approximately equal sub-meshes; the time step is $dt = 0.05$. Then follow the plots of relative distances between the smoothed reference mesh and the mesh smoothed with our method (\ref{eq:adaptedVentcell__A}--\ref{eq:adaptedVentcell__B}) for $\Delta_F$, $\Delta_A$, and $\Delta_P$. The axis $y$ shows the number of iterations and $x$ shows the relative distances (see Section \ref{sec:Evaluation} for an explanation). We want to add that for $\Delta_P$ even without decomposition the flow results in a non-smooth spiky mesh, starting at iteration 3.
}\label{fig:tri_Shrek_3721}
\end{figure}

\subsection{Evaluation}\label{sec:Evaluation}
We assess the errors introduced by the DDM in the MCF by comparing a decomposed, flown, reconstructed mesh $\tilde{M}^k$ with a reference mesh $M^k$ flown by MCF without decomposition. We compare them visually, but also in an objective manner by measuring the distances between corresponding vertices and measuring the pseudo--Hausdorff distance of the meshes $\tilde{M}^k, M^k$ with the method \verb|hausdorff_distance| from the library \verb|PyMeshLab|.

The method \verb|hausdorff_distance| returns three quantities:
\texttt{H\symbol{95}RMS} --- the root mean square of the distances of vertices of the reference mesh to the closest point (not necessarily vertex) on the candidate mesh;
\texttt{H\symbol{95}max} --- the maximum of the distances between each vertex of the reference mesh and the closest point on the candidate mesh;
\texttt{H\symbol{95}mean} --- the mean of the same distances.

We similarly denote the relative distances between corresponding vertices of the reference and candidate meshes:
\texttt{V\symbol{95}RMS} --- the root mean square of the distances of corresponding vertices;
\texttt{V\symbol{95}max} --- the maximum distance;
\texttt{V\symbol{95}mean} --- the mean of the distances.

In graphs of Figures \ref{fig:hex_F_tex_6118} and \ref{fig:tri_Shrek_3721}, we plot the relative distances, i.e., we divide the computed distances by the average edge length of the reference mesh $M^k$.

%------------------------------------------------------------------------

\begin{figure}
\centering
\includegraphics[width=0.19\linewidth]{ 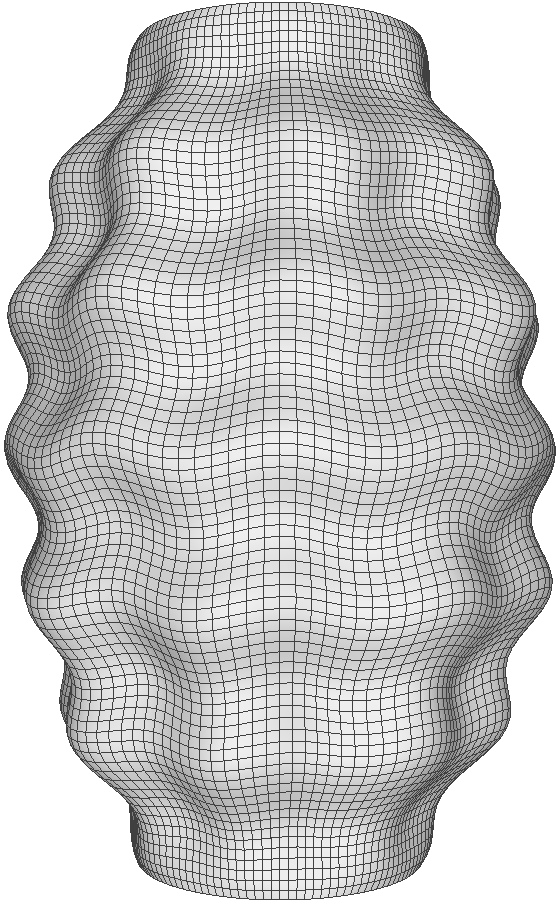}
\includegraphics[width=0.19\linewidth]{ 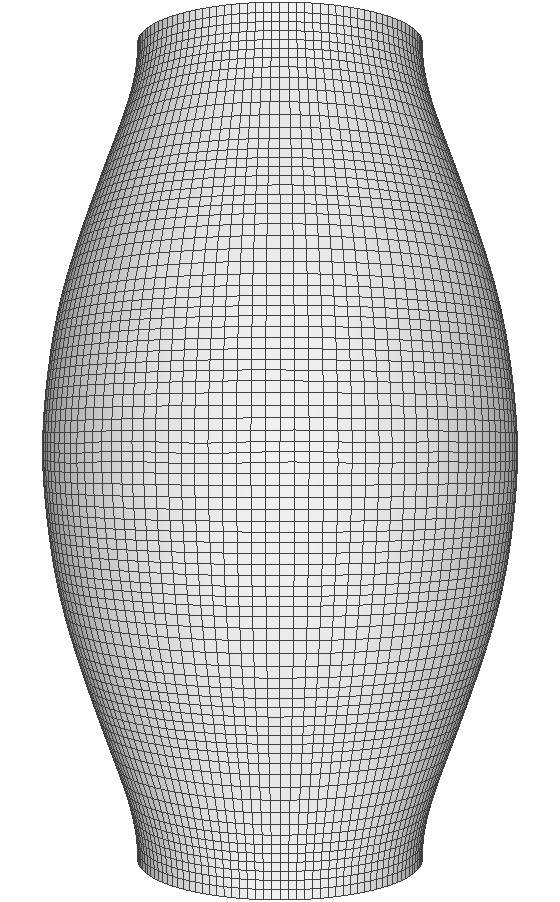}
\includegraphics[width=0.19\linewidth]{ 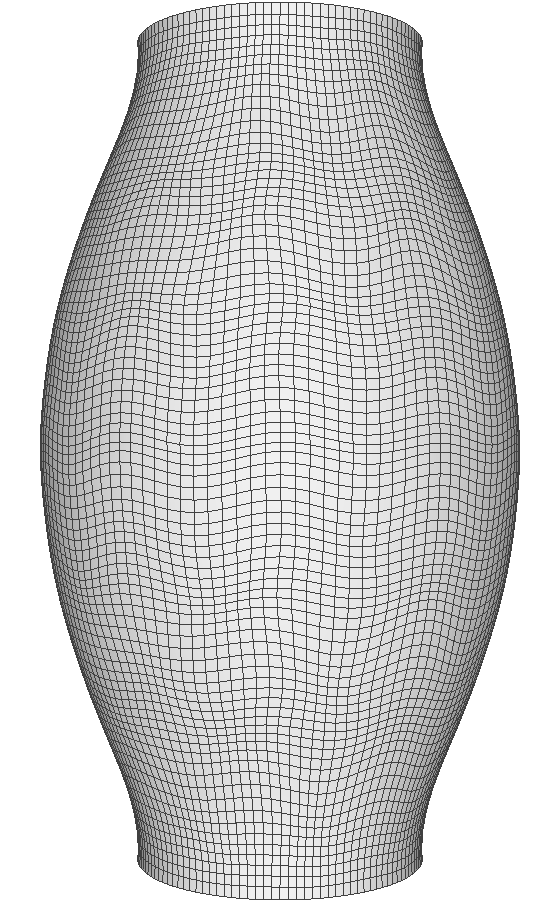}
\includegraphics[width=0.19\linewidth]{ 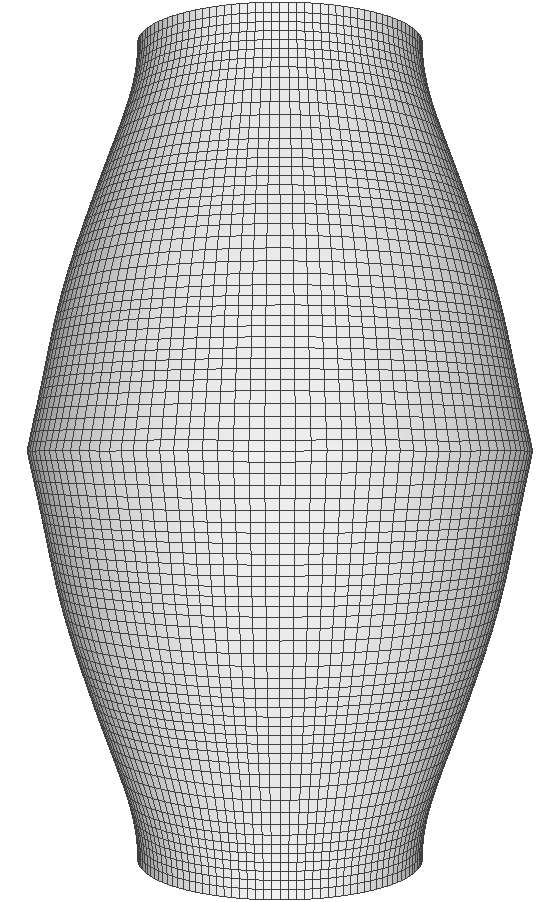}
\includegraphics[width=0.19\linewidth]{ 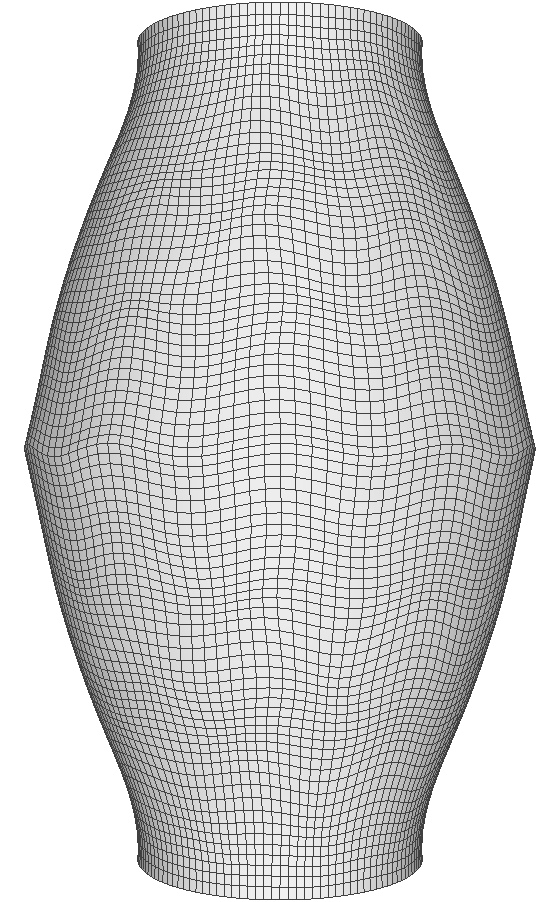}
\caption{MCF with domain decomposition of an original quad mesh (first image). We apply 5 iterations with $dt = 0.05$ and show the resulting meshes smoothed with DDM without overlap and our adapted Robin transmission conditions for $\Delta_F$ (second image) and $\Delta_P$ (third image). We also show the meshes smoothed with standard Schwarz alternating method and 1-face large overlap for $\Delta_F$ (fourth image) and $\Delta_P$ (fifth image). As can be observed, the standard Schwarz alternating method introduces significant artifacts, i.e., the vertices on the interface are moving much slower to the (eventually) minimal surface than the interior vertices.
}\label{fig:quad_cilinder_wavy_10200}
\end{figure}

%-------------------------------------------------------------------------
\section{Conclusion}
We have studied several coupling conditions together with three different Laplacians. We use three different discretization of the Laplace operator, because each one has its advantages and disadvantages. Therefore, depending on the mesh, the user may opt for one or the other. By showing, that presented domain decomposition methods work well with all the three discretization, we prove that the proposed framework has a broad application.

Concretely, the Laplacians of \cite{AlexaWardetzky2011} and \cite{PtackovaVelho2021} give rise to a mean curvature vector which has an insignificant tangential component. As a result, the MCF using these Laplacian does not deform textures mapped on the mesh, as the flow advances. The tangential shifting of vertices for Laplacian of \cite{Fujiwara1995} can be observed in Figures \ref{fig:polyMCF} and \ref{fig:quad_cilinder_wavy_10200}. This observation holds also for the MCF with domain decomposition. On the other hand, the Laplacian of \cite{Fujiwara1995} is the most stable on regular meshes.

In Section \ref{subsec:SchwarzAlternating} we present the classical Schwarz alternating method adapted to the task of MCF. We reiterate, that this method is useful, if the interface between sub-meshes is placed in an area with low mean curvature, or in area that is not expected to move significantly during the flow --- for example along a ``straight'' line, such as in the Figure \ref{fig:decomposition_poly}. On the other hand, if the interface goes through a peak, such as in Figure \ref{fig:quad_cilinder_wavy_10200}, we see a significant delay in the flow of the interface vertices in comparison to the interior vertices of the sub-meshes.
On general polygonal meshes, where it is difficult to chose a good interface between subdomains, such as in Figure \ref{fig:decomposition_poly}, the normal derivative can be difficult to define. However, the classical alternating Schwarz method can be used with good outcomes under the explained conditions on the interface.

The Ventcell transmission conditions for $q >0$ have the effect of straightening the curve of the interface between subdomains. This fact must be taken into account when choosing the proper conditions for a given mesh. Because of this drawback, this method is the one we have explored the least.

The Robin transmission conditions with normal derivative given in Definition \ref{def:newNormalDerivative} and $p^*$ given in (\ref{eq:ourP}) guarantee that the overall shape of evolving surface mesh will mimic the shape of corresponding family of meshes when MCF without decomposition is applied. The interface does not have to be a straight curve, however it must allow for a well defined normal derivative. Moreover, we have shown, that this method allows for MCF with decomposition without introducing artifacts such as texture deformation. Even more, our adaptation is also reliable when the initial mesh has different levels of details in certain areas, as illustrated in Figure \ref{fig:teaser}.

To conclude, we believe that our exposition convinced the interested reader, that domain decomposition methods are worth exploring in the setting of mesh smoothing. Principally, because in the age of parallel computing, they offer a guaranteed speed--up of methods that involve the Poisson equation, as explained in Section \ref{sec:Implementation}.

%%%%%%%%%%%%%%%%%%%%%%%%%%%%%%%%%%%%%%%%%%%%%%%%%%%%%%%%%%%%%%%%%%%%%
%%%%%%%%%%%%%%%%%%%%%%%%%%%%%%%%%%%%%%%%%%%%%%%%%%%%%%%%%%%%%%%%%%%%%
%\section*{References} % needed on some systems
\bibliography{GMP2026}

\end{document}